\documentclass[11pt]{amsart}
\usepackage{geometry}                
\usepackage{graphicx}
\usepackage{amssymb}
\usepackage{epstopdf}
\newcommand{\eps}{\varepsilon}
\newcommand{\diag}{\mbox{diag}}
\newcommand{\tr}{\mbox{tr}}
\newcommand{\T}{\mathbb{T}}
\newcommand{\R}{\mathbb{R}}
\newcommand{\Z}{\mathbb{Z}}

\title{Interaction of two systems with saddle-node bifurcations on invariant circles\\ I. Foundations and the mutualistic case.}
\author{C.Baesens \and R.S.MacKay}
\address{Mathematics Institute, University of Warwick, Coventry CV4 7AL, U.K.}

\begin{document}
\maketitle
\begin{abstract}
The saddle-node bifurcation on an invariant circle (SNIC) is one of the codimension-one routes to creation or destruction of a periodic orbit in a continuous-time dynamical system.  It governs the transition from resting behaviour to periodic spiking in many class I neurons, for example.  Here, as a first step towards theory of networks of such units the effect of weak coupling between two systems with a SNIC is analysed.  Two crucial parameters of the coupling are identified, which we call $\delta_1$ and $\delta_2$.  Global bifurcation diagrams are obtained here for the ``mutualistic''  case  $\delta_1\delta_2>0$.  According to the parameter regime, there may coexist resting and periodic attractors, and there can be quasiperiodic attractors of torus or cantorus type, making the behaviour of even such a simple system quite non-trivial.  In a second paper we will analyse the mixed case $\delta_1\delta_2<0$ and summarise the conclusions of this study.
\end{abstract}
\noindent{\bf Keywords}:
saddle-node on invariant circle; bifurcation; class I neuron, Josephson junction.\\
\noindent{\bf MSC}: 37Exx, 37N25, 37Gxx

\section{Introduction}   

The saddle-node on an invariant circle (SNIC) bifurcation is one of the basic scenarios for creation of a periodic orbit in smooth continuous-time dynamical systems (number three in the list of Andronov and Pontryagin ~\cite{AP}). It goes under various other names, such as saddle-node on a limit cycle (SNLC), homoclinic to a saddle-node, or saddle-node loop (though we reserve the latter name for a codimension-2 case, e.g.~$Z$ points of \cite{BGKM}).

The SNIC bifurcation is for example the scenario by which many class I neurons are believed to make the transition from resting behaviour to periodic spiking ~\cite{EK,RE} (this was also proposed by RSM to the neurophysiologist H.B.Barlow in 1986).  It underlies some regenerative excitable chemical systems.  It occurs in mechanical systems  too, like the damped pendulum with torque for sufficiently strong damping, or its Josephson junction analogue ~\cite{LHM}.

Mathematically, a SNIC for a $C^r$ ($r\ge 2$) vector field $\dot{X}=V(X)$ on a manifold is an elementary saddle-node (an equilibrium with a simple eigenvalue $0$, no other eigenvalues on the imaginary axis and non-zero quadratic coefficient in the null direction) with a trajectory that is asymptotic to the saddle-node in both directions of time (``homoclinic") along the null direction; it follows that the resulting invariant circle is $C^r$.  Thus, taking a coordinate $x$ along the circle, with $0$ representing the saddle-node equilibrium, and suitable orientation, the vector field $v$ on the circle satisfies 
\begin{equation*}
\dot{x} = v(x) = Cx^2 + o(x^2)
\end{equation*}
 as $x\to0$ for some $C>0$ and $v(x)>0$ for $x\ne0$.   Without loss of generality, by scaling $x$ we take $C=1$.  We denote the length of the circle (in the rescaled coordinate $x$) by $L$ (an alternative would be to scale the length of the circle to 1 and keep a general value of $C$).

Consider $C^r$ perturbation of a SNIC by a parameter $\mu$, with $\mu=0$ representing the unperturbed case.  By the elementarity condition, the invariant circle is normally hyperbolic and so persists for all $C^r$-small perturbations. A parameter-dependent extension of the coordinate $x$ can be chosen so that the length of the circle remains $L$ and the vector field on it is a $C^r$-small perturbation of that for $\mu=0$.  By a parameter-dependent shift of origin to remove the linear term (using the implicit function theorem), the perturbed vector field has the local form 
\begin{equation}
\dot{x} = a(\mu) + b(\mu) x^2 + o(x^2)
\label{eq:sne}
\end{equation} 
for some smooth functions $a, b$ of $\mu$, with $a(0)=0, b(0)=1$.  

It is convenient for later purposes to make $b(\mu)$ precisely $1$ for $\mu$ small, by a coordinate change $X(x, \mu)$ preserving the length $L$.  This can be achieved as in Appendix A.
We will suppose that $a'(0)\ne 0$ and thus for small $\mu$ we can use $a$ as parameter instead of $\mu$, so without loss of generality we have 
\begin{equation}
\dot{x} = \mu+x^2 +o(x^2)\label{eq:dotxmux2}
\end{equation}  
as $x\to0$.  For $\mu<0$ it has two equilibria: a sink and a source (if the normal directions to the circle are attracting, the case of most relevance, then in the full state space these are a sink and a saddle).  For $\mu>0$ the circle is a periodic orbit whose period $T(\mu)$ is asymptotic to $\pi/\sqrt{\mu}$ as $\mu \to 0$, which spends all but a bounded amount of its period in any neighbourhood of $0$.

An explicit example of a family with a SNIC bifurcation is 
$$\dot{\phi} = \mu + 2(1-\cos{\phi})$$ 
on a circle of length $2\pi$.  Another artificial-looking but useful example is 
$$\dot{x} = \mu + x^2$$ 
on the real line union a point at infinity, interpreted as the projective real line (set of lines $\ell$ through the origin in a plane); the projective real line is diffeomorphic to a circle with the coordinate $x$ representing the slope of the line $\ell$; one could write $x = \tan{\theta}$ with $\theta$ considered modulo $\pi$ (half a revolution brings $\ell$ to itself) and then 
$$\dot{\theta} = \mu \cos^2{\theta} + \sin^2{\theta}$$
on a circle of length $\pi$, which is similar to the first example (put $\phi=2\theta$), cf.~\cite{Iz}.

The goal of this paper is to find what happens when two SNIC bifurcations are coupled.  This is step one towards finding how a network of class I neurons or Josephson junctions behaves, as addressed for example in chapter 8 of \cite{HI}.  Although it sounds a simple problem, we have not found a rigorous treatment in the literature and our analysis and its solution are remarkably complicated.  An example was treated numerically by \cite{GK}.  We treat the general case from the point of view of determining the minimal structure that the bifurcation diagrams must possess.  

We identify two crucial coupling parameters $\delta_1$, $\delta_2$.  We distinguish two principal cases: ``mutualistic'', with $\delta_1\delta_2>0$, and ``mixed'', with  $\delta_1\delta_2<0$.  The mutualistic case can be further decomposed into ``mutually excitatory'': $\delta_1,\delta_2>0$, and ``mutually inhibitory'': $\delta_1,\delta_2<0$.  We analyse the simplest mutualistic cases in this paper, but even they turn out to be quite complicated.  In a companion paper, Paper II, we will treat the mixed case and summarise the conclusions of the study.

We will assume the uncoupled systems are $C^r$ with $r$ large enough that every operation we will perform makes sense ($r=5$ suffices). The coupling we treat is of the form of a $C^r$-small perturbation to the vector field given by the product of two SNIC bifurcations.  This could represent the effect of gap-junction coupling (``electrical synapse'')  between neurons (for an introduction to neuroscience see \cite{NMW}).  We will show in Section~\ref{sec:insideTriangle} that gap-junction coupling produces the mutually excitatory case.  On the other hand, chemical synaptic coupling requires further analysis, because it may involve adding additional degrees of freedom to represent the effects of the neuro-transmitters and time delays to represent signal travel times.  The first consideration does not make much difference, because near the bifurcation the timescale for the electric dynamics is longer than for relaxation of the channels in response to the neuro-transmitters so the latter can simply be added into the normally hyperbolic directions.  The second makes the dynamics infinite-dimensional and although it is likely that again the effects can simply be added into the normally hyperbolic directions, we did not pursue this yet.  Similarly, the coupling can model electrical coupling between Josephson junctions.  For an introduction to the theory and experiments on Josephson junction arrays, see ~\cite{M,U} respectively.

Our analysis uses heavily terminology and results from section 4 of \cite{BGKM} (on bifurcations for flows on a 2-torus) and some from \cite{BM} (on coupling of a saddle-node periodic orbit with an oscillator).  In particular, we recall some key concepts right now.
A {\em Poincar\'e flow} on $\T^2$ is one with a {\em global cross-section}, i.e.~a transverse section such that every forward and backward orbit crosses it.  With respect to a choice of coordinates $(x_1,x_2)$ on the universal cover of the torus (i.e.~consider $\T^2 = \R^2/(L_1\Z\times L_2\Z)$, where $L_j$ are the lengths of the cycles in the coordinate directions), the {\em homology direction} of a forward orbit of a flow on $\T^2$ is the limit of the unit vector in the direction of the vector $V$ of (signed) numbers of revolutions in $x_1$ and $x_2$ as time goes to $+ \infty$ (or $0$ if $V$ does not go to infinity).  The {\em winding ratio} is the homology direction modulo reflection through $0$.  For a Poincar\'e flow, every orbit has the same homology direction and it is non-zero.  We denote Poincar\'e flows by $P$.
A {\em Cherry flow} is one with a homotopically non-trivial transverse section $\Sigma$ and a direction of time such that the orbits of a non-empty subset $\Sigma'$ return to $\Sigma$ under the flow, the induced map $g: \Sigma'\to\Sigma$ is continuous, and $\lim_{x\to l}g(x)=\lim_{x\to r}g(x)$ for all gaps $[l,r]$ (components of $\Sigma\backslash\Sigma'$).
 Every unbounded orbit of a Cherry flow has the same non-zero homology direction.  We denote Cherry flows by $C$ (in \cite{BGKM}, $C$ denotes a larger class).
If the homology direction of a flow is that of an integer vector $(p_1,p_2)$ with no common factors, we say the flow is {\em partially mode-locked of type} $(p_1,p_2)$.
A flow is {\em fully mode-locked} if every orbit has homology direction $0$ (equivalently if every orbit is bounded on the universal cover).  We denote fully mode-locked flows by $F$ (or $FML$).  A non-contractible closed curve on the torus is called {\em rotational}.

\section{Product system}    
\label{sec:prod}

Denote the coordinates of the two systems by $x_j$, $j=1,2$, with lengths $L_j$ for one revolution, their parameters by $\mu_j$ and their (uncoupled) vector fields by 
$$\dot{x}_j = v_j(x_j,\mu_j) = \mu_j + x_j^2 + O(x_j^3)$$
with $v_j(x_j,0)>0$ for all $x_j\ne0$.   We interpret the remainder term in at least $C^1$, thus in particular there are $K>0$, $M>0$ such that $|v_j-\mu_j-x_j^2|\le K|x_j|^3$, $|v_j'-2x_j|\le 3K|x_j|^2$, for $|x_j|\le M$.

It will be convenient to suppose that $\frac{\partial v_j}{\partial \mu_j} \ge c> 0$ for all $x_j$, not just for $x_j$ near $0$ (say $c=\frac12$).  This can be achieved by parameter-dependent coordinate changes $X_j(x_j,\mu_j)$ as in Appendix B.

\begin{figure}[htbp]     
   \centering
   \includegraphics[width=3in]{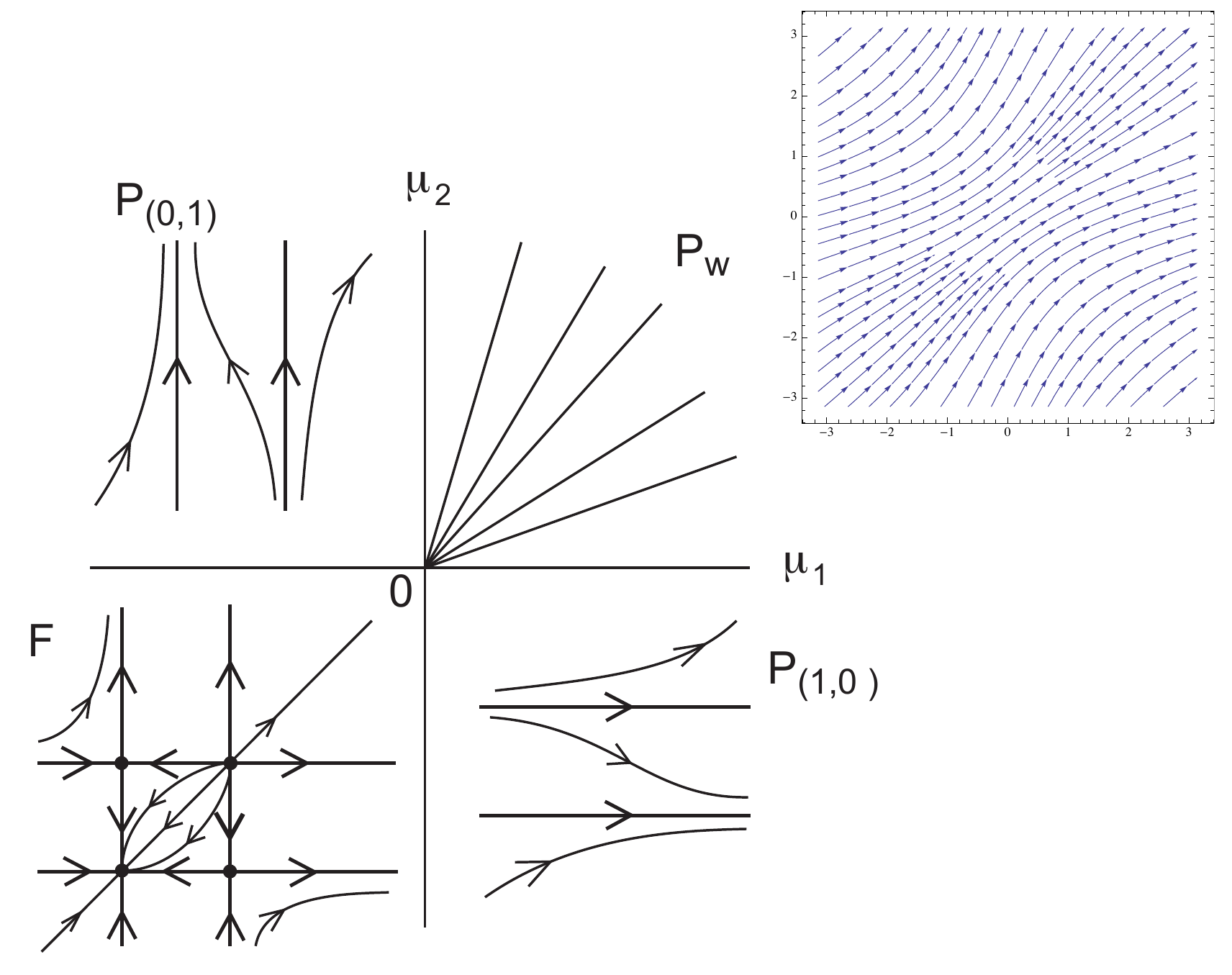}
    \caption{Bifurcation diagram in the parameter plane $(\mu_1,\mu_2)$ for two uncoupled SNICs, with global phase portraits.}
   \label{fig:bifprod}
\end{figure}

The product of two circles is a 2-torus.  We obtain the bifurcation diagram of Figure~\ref{fig:bifprod} for the product system in the plane of $\mu = (\mu_1,\mu_2)$, with global phase portraits as indicated.  
In the positive quadrant, the flow is smoothly conjugate to that of a constant vector field on the unit torus which varies smoothly with parameters.  For the conjugacy one can use the fractions $\tau_j(x_j)$ of the period $T_j$ traversed from a reference point (say $x_j=0$).  The constant vector field is $(1/T_1,1/T_2)$, and it has asymptotic expression $(\sqrt{\mu_1},\sqrt{\mu_2})/\pi$ as $\mu\to 0$.  In particular, the flows in the positive quadrant are Poincar\'e flows  with homology direction $(T_2,T_1)/\sqrt{T_1^2+T_2^2}$ which varies smoothly and at non-zero rate with the ratio $\mu_1:\mu_2$.
In the negative quadrant, the flow is fully mode-locked, with four invariant circles ($x_j \approx \pm \sqrt{-\mu_j}$) intersecting pairwise at four equilibria.  We call phase portraits topologically equivalent to this, {\em basic tartan}.  The boundaries of the negative quadrant constitute two simultaneous saddle-node bifurcations of equilibria (sne), which merge into a degenerate one at $\mu=(0,0)$.  
In the bottom right quadrant, the flow is a Poincar\'e flow
of type $(1,0)$, with repelling and attracting periodic orbits at $x_2 \approx \pm \sqrt{-\mu_j}$.  Similarly, in the top left quadrant, we have Poincar\'e flow of type $(0,1)$.
The boundaries of the positive quadrant (minus the vertex $0$) correspond to (elementary) saddle-node periodic orbits (snp) (periodic orbits with a Floquet multiplier $+1$).

\section{Effect of weak coupling}   

\subsection{First steps}   
\label{sec:firststeps}

The invariant 2-torus of the uncoupled system is normally hyperbolic for $\mu$ small enough, so persists under small smooth perturbation \cite{F}, and the vector field on it is a small smooth perturbation of the uncoupled case, in general $\mu$-dependent.  We restrict attention from now on to a neighbourhood of $\mu=0$ where the above holds.  Let us denote the perturbation size in $C^r$, for some $r\ge2$, by $\delta$.  In particular this implies that the changes to $v_1$ and $v_2$ and to their first and second derivatives are at most $\delta$ (actually  in Appendix E we will require a bound on the third derivative of the perturbation, which does not need to be as strong as $\delta$ but for convenience we will assume that too).

The perturbed system is 
\begin{align*}
\dot{x}_1=\tilde{v}_1(x_1,x_2)&=v_1(x_1)+\mathcal{O}(\delta)&&\text{with}&v_1(x_1)&=\mu_1+x_1^2+\mathcal{O}(x_1^3)\\
\dot{x}_2=\tilde{v}_2(x_1,x_2)&=v_2(x_2)+\mathcal{O}(\delta)&&\text{with}&v_2(x_2)&=\mu_2+x_2^2+\mathcal{O}(x_2^3)\,.
\end{align*}
\begin{figure}[htbp]       
   \centering
   \includegraphics[height=2in]{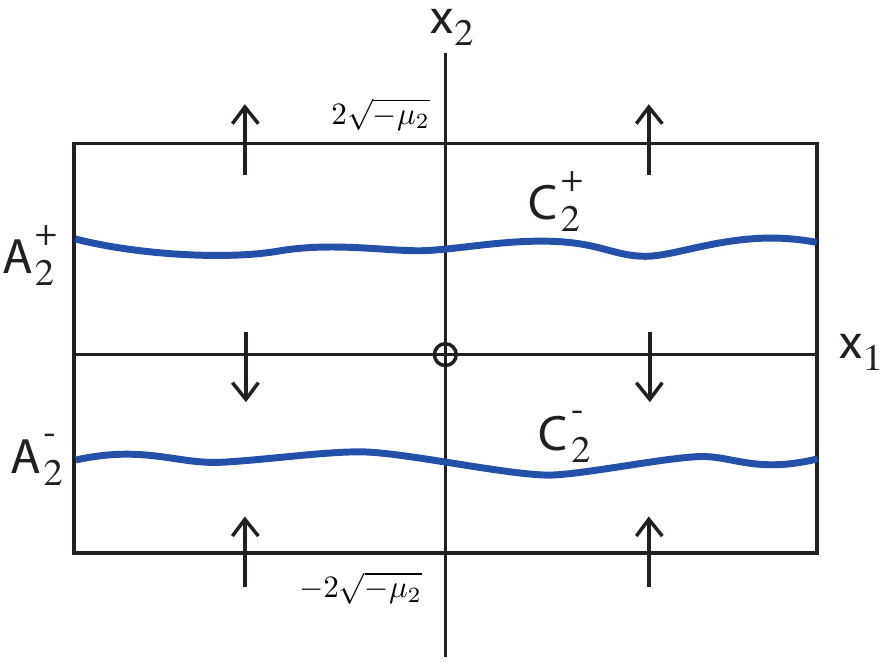} 
   \caption{Annuli $A^+_2:=\{0\le x_2\le2\sqrt{-\mu_2}\}$ and $A^-_2:=\{-2\sqrt{-\mu_2}\le x_2\le0\}$, and resulting homotopically non-trivial repelling and attracting sets $C^\pm_2$ (drawn as $C^1$ invariant circles here but could take other forms as in Figure~\ref{fig:figureX}).}
   \label{fig:annuli}
\end{figure}
For $-\frac{1}{16K^2}<\mu_2<-\delta$ we have $\dot{x}_2<0$ on $x_2=0$ and $\dot{x}_2>0$ for $|x_2|\ge2\sqrt{-\mu_2}$.
Let $A_2^\pm$ be the annuli as shown in Figure~\ref{fig:annuli}.  Then, defining $\phi$ to be the flow of the above differential system, $C_2^-:=\cap_{t>0}\phi_tA_2^-$ and $C_2^+:=\cap_{t<0}\phi_tA_2^+$ are homotopically nontrivial attracting, respectively repelling sets. 

\begin{figure}[htbp] 
   \centering
   \includegraphics[width=5.5in]{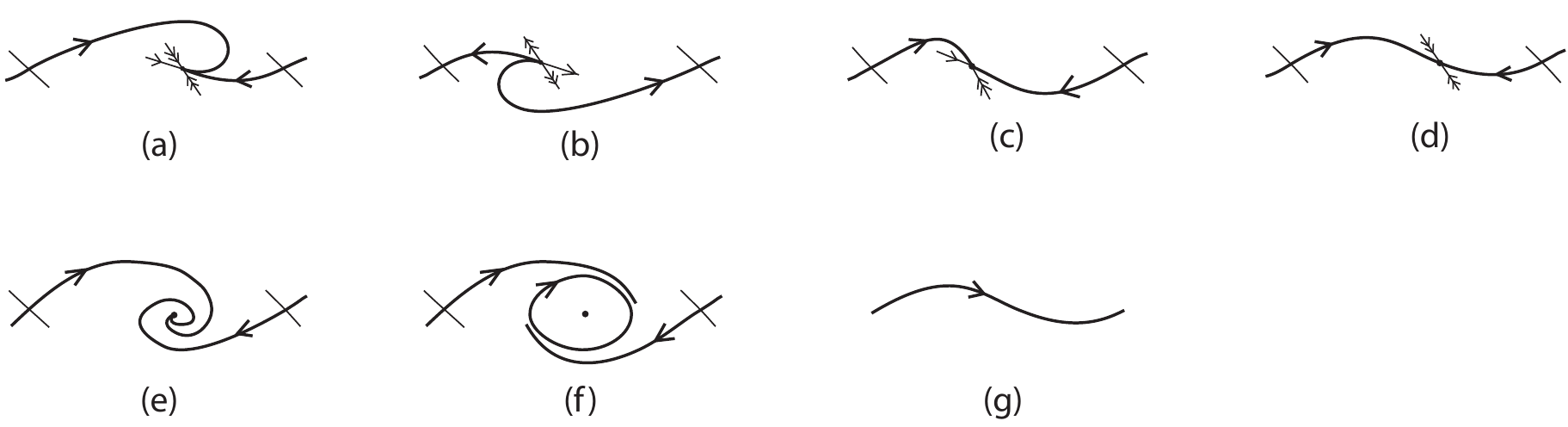} 
   \caption{Some possible forms for the maximal invariant sets in the annuli. Only (d) and (g)  are $C^1$. Case (f) is not a $C^0$ circle.} 
   \label{fig:figureX}
\end{figure}

If we neglect the perturbation, $C^\pm_2$ are just the circles $x_2\approx\pm\sqrt{-\mu_2}$ given by the zeroes of $v_2$. 
Under perturbation, we will find regions in which they persist to $C^1$ invariant circles (either periodic orbits or chains of connecting orbits between equilibria), regions in which they are $C^0$ invariant circles connecting equilibria but not necessarily $C^1$, and in Paper II regions in which they are not even $C^0$ circles.  Figure~\ref{fig:figureX} shows some possible forms for $C^\pm_2$. 

We begin by taking $\mu_2\le-2\delta$ and showing in Appendix C that the $x_2$-nullcline consists of two $C^1$ graphs over $x_1$, lying within $\frac23\sqrt{-\mu_2}\le\pm x_2\le2\sqrt{-\mu_2}$ and having small slope (at most $\delta/\sqrt{-\mu_2}$).  The analogous result holds for the $x_1$ nullcline when $\mu_1\le-2\delta$.

\begin{figure}[htbp] 
   \centering
   \includegraphics[width=5.5in]{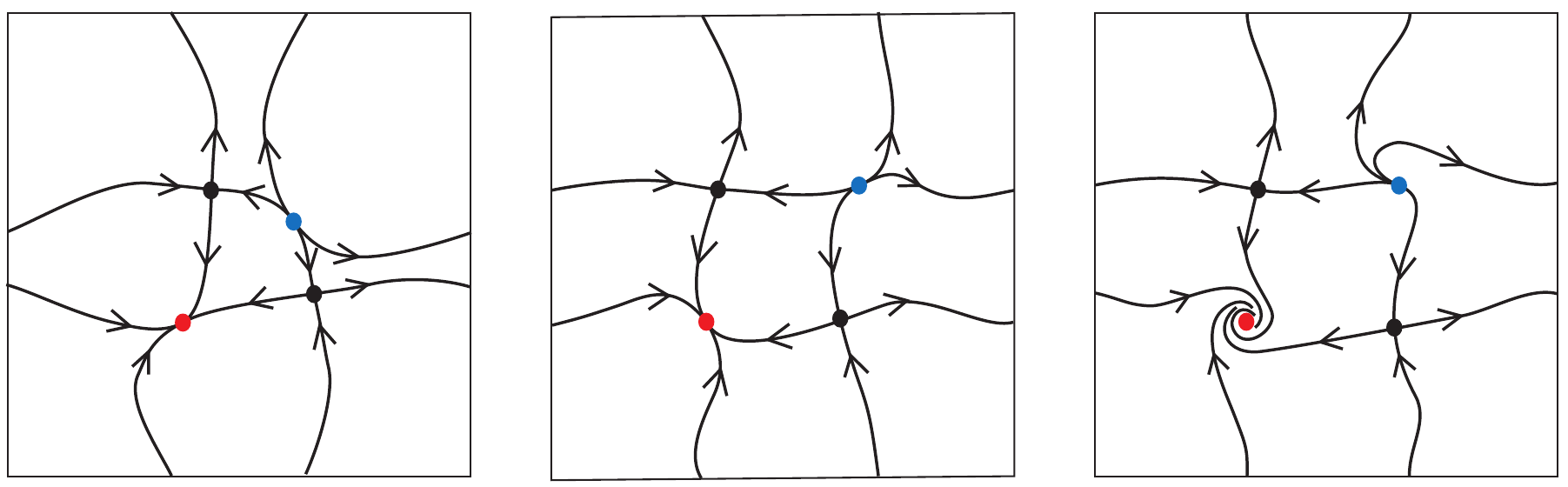} 
   \caption{Some possible realisations of basic tartan in $C^1$. The third one has $C^\pm_2$ and $C^-_1$ non-$C^1$ circles.}
   \label{fig:figureY}
\end{figure}

In particular, in the region where both $\mu_j\le-2\delta$, the nullclines intersect in precisely four points.  As the derivative of $\tilde{v}$ is close to $\text{diag}\,(2x_1,2x_2)$, they are two saddles, a source and a sink, arranged just as for the unperturbed case.  We leave out the detailed justification.  The invariant manifolds of the saddles leave close to horizontal and vertical and because of the signs of the components of $\tilde{v}$ between the nullclines, they are obliged to fall into the source or sink in topologically the same way as for the unperturbed case.
Thus when both $\mu_j\le-2\delta$, we continue to obtain a basic tartan phase portrait  at the $C^0$ level, though at the $C^1$ level it can take forms like those in Figure~\ref{fig:figureY}.  In particular, $C^\pm_2$, $C^\pm_1$ are all at least $C^0$ circles , but not necessarily $C^1$ because of the ways the saddle manifolds may meet at the source or sink.

Next we analyse $\mu_1\ge-2\delta$, $\mu_2\ge-C\delta$ for some $C$ greater then 2 and prove that $C^\pm_2$ are $C^1$ invariant circles there, using normal hyperbolicity theory.
The normal linearised dynamics for the unperturbed case
$$\dot{\delta x_2}=v_2'(x_2)\,\delta x_2\approx\pm2\sqrt{-\mu_2}\,\delta x_2$$
is hyperbolic (repelling for $+$, attracting for $-$).  The character of the tangential linearised dynamics
$$\dot{\delta x_1}=v_1'(x_1)\,\delta x_1$$
depends on the sign of $\mu_1$.

For $\mu_1\ge0$, although tangent orbits grow a lot for $0<x_1\ll L_1/2$, the growth is all cancelled out by contraction for $0<L_1-x_1\ll L_1/2$, producing zero Lyapunov exponent.  Thus taking Fenichel's approach \cite{F} to normal hyperbolicity theory, time-averaged tangential contraction or expansion rates are less than normal ones and so the circles persist to nearby $C^1$ invariant circles on adding $C^1$ small enough perturbation.

For $\mu_1<0$, there are two equilibria $x_1\approx\pm\sqrt{-\mu_1}$ on the invariant circles of the unperturbed system, with Lyapunov exponents approximately $\pm 2\sqrt{-\mu_1}$. All other orbits are heteroclinic to these so have forward Lyapunov exponent $-2\sqrt{-\mu_1}$ and backward Lyapunov exponent $+2\sqrt{-\mu_1}$.  Thus tangential contraction or expansion rates are weaker than the normal ones if $\mu_2<\mu_1<0$ and so under this condition the circles persist to nearby $C^1$ invariant circles on adding $C^1$ small enough perturbation for parameters in this region.

\begin{figure}[htbp] 
   \centering
   \includegraphics[width=2.5in]{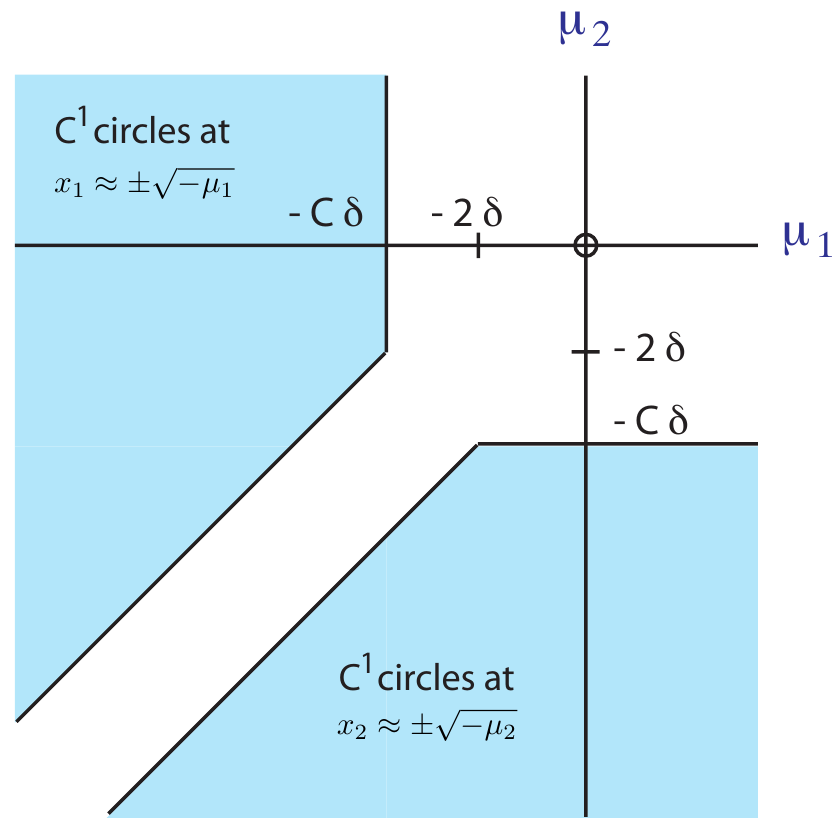} 
   \caption{There are $C^1$ rotational invariant circles in the shaded areas.}
   \label{fig:figureW}
\end{figure}

To quantify what counts as small enough perturbation, however, is not so easy.  In Appendix D we show there is $C>2$ such that the $C^1$ invariant circles persist in at least the region $\mu_2\le-C\delta$, $\mu_1\ge\mu_2+(C-2)\delta$, sketched in Figure~\ref{fig:figureW}. 

A $C^1$ invariant circle may either be a periodic orbit or a $C^1$ string of connections between equilibria.  The case of periodic orbit happens for $\mu_1>C\delta$ (with $\mu_2<-C\delta$) because then $\dot{x}_1>0$.  The case of a $C^1$ string of connections between equilibria happens for $\mu_1\le-2\delta$  because then there are two equilibria in each annulus.
 
Outside the region where normal hyperbolicity theory applies, however, it could happen that the attracting or repelling sets have the forms of Figure~\ref{fig:figureX}(a), (b) respectively, or more complicated, e.g.~Figure~\ref{fig:figureX}(e),(f).  Transition from $C^1$ circle (Figure~\ref{fig:figureX}(d)) to other phase portraits of Figure~\ref{fig:figureX} would proceed via use of the fast direction to the sink (respectively source) (Figure~\ref{fig:figureX}(c)) or via node to focus transition leading to Figure~\ref{fig:figureX}(e) and possible subsequent Hopf bifurcation leading to Figure~\ref{fig:figureX}(f).  For an example of such transitions see \cite{BM}.

Similar remarks hold for $\mu_1<-C\delta$ and vertical invariant circles (see Figure~\ref{fig:figureW}).

For $\mu_1+\mu_2 > C \delta$ there is a global cross-section, e.g.~$\frac{x_1}{L_1}+\frac{x_2}{L_2} = 0$, because $\tilde{v}_1+\tilde{v}_2>0$ everywhere in that region of parameter space.  Thus we obtain Poincar\'e flows.  The homology direction goes continuously from $(1,0)$ to $(0,1)$ from the lower right to the upper left, but generically locks to rational values.  In fact it changes monotonically with a parameter $\lambda$ along the lines $\mu_1 = \frac{K}{2}-\lambda, \mu_2 = \frac{K}{2}+\lambda$, $K$ constant ($K>C \delta$), because for the unperturbed system the component of the derivative $\frac{\partial v}{\partial \lambda}$ in the direction $v^\perp = (-v_2,v_1)/|v|$ (using Euclidean norm $|v| = \sqrt{v_1^2+v_2^2}$) is at least $c(v_1+v_2)/|v|$ (where $c$ is as in Section~\ref{sec:prod}), which is positive in this region, and small perturbation can not change its sign, so increasing $\lambda$ turns the vector field in the positive (anticlockwise) direction.  On the boundaries of the regions of Poincar\'e flow of rational type $(p,q)$ there is a saddle-node periodic orbit.

\begin{figure}[htbp] 
   \centering
   \includegraphics[width=3.5in]{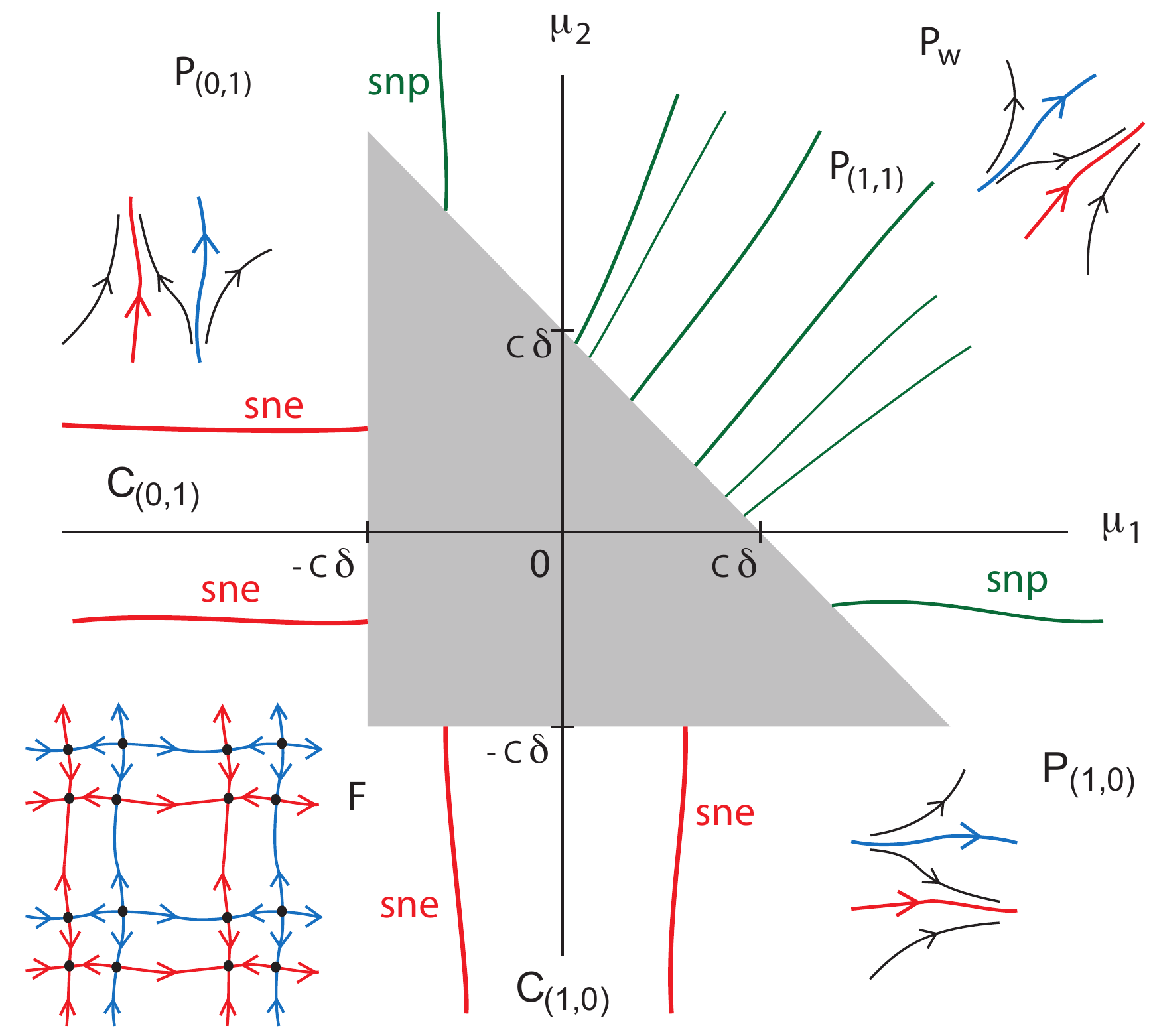} 
   \caption{Bifurcation diagram of coupled system outside a triangle in parameter space. The phase portrait in the fully mode-locked region $F$ is $C^0$-correct but generically takes $C^1$ forms like those of Figure~\ref{fig:figureY}.  Various phase portraits are possible in the regions of Cherry flow; see Figure~\ref{fig:insidetriangle1}, for example.}
   \label{fig:outside}
\end{figure}

Putting the results of the previous two paragraphs together, we have a region of full mode-locking which includes $\mu_1, \mu_2 < -C\delta$, we have a region of Poincar\'e flow of type $(1,0)$ which includes $\mu_2<-C\delta, \mu_1+\mu_2>C\delta$, and one of type $(0,1)$ which includes $\mu_1<-C\delta, \mu_1+\mu_2>C\delta$, and a region of Poincar\'e flows with homology direction varying monotonically for $\mu_1+\mu_2>C\delta$.  See Figure~\ref{fig:outside}.

For $\mu_2<-C\delta$ there are precisely two curves of saddle-node equilibria, graphs over $\mu_2$, separating $\mu_1<-C\delta$ from $\mu_1+\mu_2>C\delta$ and each creating a pair of equilibria on one of the horizontal invariant circles.  This is because the equations for a saddle-node equilibrium are $\tilde{v}_1=0$, $\tilde{v}_2=0$ and $\det D\tilde{v} = 0$, which can be written approximately as $\mu_1 \approx -x_1^2$, $x_2^2 \approx -\mu_2$, $4x_1x_2 \approx 0$; using the implicit function theorem, the second equation determines $x_2$ as either of two functions of $(x_1, \mu_1, \mu_2)$ near $x_2 = \pm \sqrt{-\mu_2}$ for $\mu_2<-C\delta$; substituting these into the third determines $x_1 \approx 0$ as one of two functions of $(\mu_1, \mu_2)$; substituting for $x_2(x_1,\mu_1,\mu_2)$ and then $x_1(\mu_1,\mu_2)$ into the first provides $\mu_1 \approx 0$ as either of two functions of $\mu_2$.
Without further hypotheses it is not possible to say which saddle-node bifurcation happens first (indeed the curves could cross), but in between them the flow is Cherry flow of type $(1,0)$.  Similarly, there are two curves of saddle-node equilibrium in $\mu_1<-C\delta$, graphs over $\mu_1$, in between which we have Cherry flow of type $(0,1)$.

The phase portrait in the region $\mu_j \le -C\delta$, $j=1,2$, is a basic tartan as we already discussed.  It is not guaranteed to remain like this in the whole of the full mode-locked region, however.  For example, if the two sne curves in $\mu_2 \le -C\delta$ cross then it could easily happen that to the left of this there is a heteroclinic bifurcation $D \to A_{01}$ (see Fig.~\ref{fig:ABCDlabels} for the notation), which would give rise to a skewed tartan (two invariant horizontal circles and two invariant circles of type $(1,1)$).  Nevertheless, under generic hypotheses to be formulated at the end of Sec.~3.2 ($\delta_1,\delta_2 \ne 0$), we will prove in Sec.~4.3 that the phase portrait is a basic tartan in all the part of the full mode-locked region with $\min_{j\in\{1,2\}} |\mu_j| \gg \delta^2 (\log{1/\delta})^4$ for $\delta_1\delta_2>0$ (the case $\delta_1\delta_2<0$ will be addressed in Paper II).

Thus we have determined the principal features of the perturbed bifurcation diagram outside the triangle in parameter space shown in Figure~\ref{fig:outside}.

\subsection{Inside the triangle} 
\label{sec:insideTriangle} 

\begin{figure}[htbp] 
   \centering
   \includegraphics[width=2.5in]{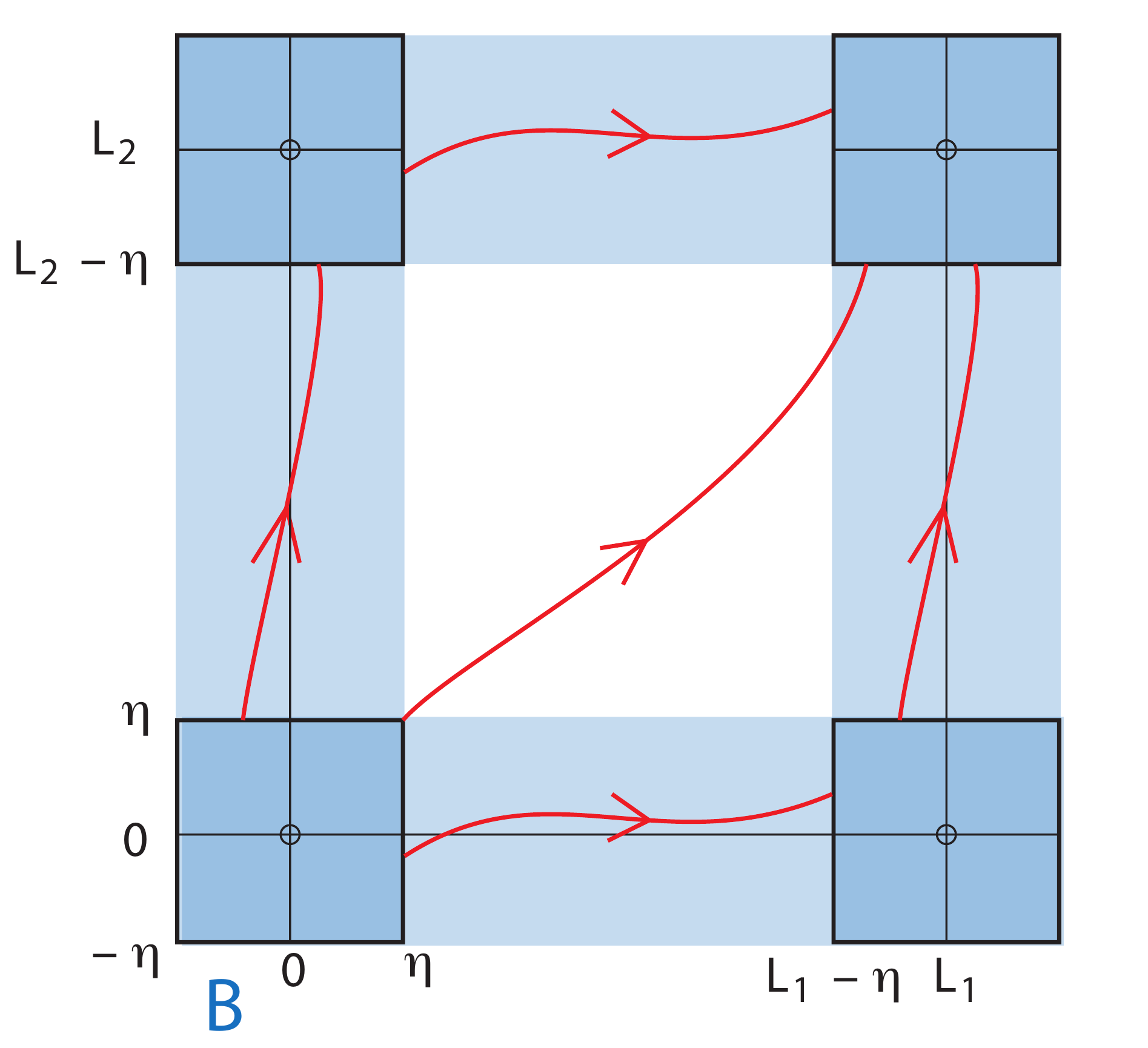} 
   \caption{Strips and box $B$.}
   \label{fig:box&strips}
\end{figure}

To analyse what happens inside the triangle in parameter space of Figure~\ref{fig:outside}, we divide the torus into the strips $|x_1|\le\eta,|x_2| \le \eta$ for some $\eta$ small, their intersection $B$ and the complement (Figure~\ref{fig:box&strips}).

The idea is that $\eta$ should be small enough so we can accurately use second order Taylor expansion of $v_j$ about $x_j=0$ inside the strip $|x_j|\le \eta$, yet considerably larger than $\sqrt{C\delta}$ so that the effects of $\mu = (\mu_1,\mu_2)$ and $\delta$ outside the strips are relatively small, in particular so that $\tilde{v}_j >0$ outside $|x_j|\le \eta$ and there are no equilibria outside the intersection box.  We will start with $\eta = k \delta^{1/3}$ for a small constant $k$, which is the largest that suffices for neglect of the cubic and higher terms compared with $\mu$ when $\mu$ is inside the triangle, but we will at later stages reduce or increase it, so we leave $\eta$ explicitly in most formulae.  

First we consider the vector field $\tilde{v}$ inside the intersection box $B: |x_j|\le \eta$ for both $j=1,2$.  It is a small perturbation of $\dot{x}_1 = \mu_1+x_1^2, \dot{x}_2 = \mu_2 + x_2^2$.  Let us Taylor expand about $x_1,x_2=0$.  

The quadratic terms $q(x)$ of $\tilde{v}$ can be transformed to $X_1^2 + \eps_1 X_2^2, X_2^2+\eps_2 X_1^2$ by putting $x = MX$ with $M$ some near-identity matrix, depending on $\mu$, determined by eliminating the $X_1X_2$ terms in the quadratic part $M^{-1}q(x)$ of $\dot{X}$ and making its ``diagonal" coefficients be $1$.  The idea is that the derivative of the map from the matrix elements $M_{11},M_{12},M_{21},M_{22}$ to the coefficients of the $X_1^2$ and $X_1X_2$ terms in $\dot{X}_1$ and the $X_1X_2$ and $X_2^2$ terms in $\dot{X}_2$ is near the matrix $\diag(1,2,2,1)$ which is invertible.  Thus $M$ comes out within $\delta$ of the identity (since $\delta$ bounds the size of the perturbation to the quadratic terms in $\dot{x}$) and the coefficients $\eps_j$ above are also of order $\delta$.  It is not strictly necessary to have reduced the quadratic terms in this way, since the quadratic part of the perturbation would have turned out to be negligible anyway, but it is tidier to eliminate terms when one can rather than have to bound their effects.

The above reduction of the quadratic terms can be equally well achieved with a global coordinate change on the torus by using instead 
$$
x = X + (M-I)\left[\begin{array}{c}
k_1^{-1} \sin{k_1X_1} \\
k_2^{-1} \sin{k_2X_2}
\end{array}\right],
$$
with $k_j = 2\pi/L_j$, which is within order $\delta$ of the identity.  Conjugating by this changes the vector field by order $\delta$, which is the same size as the perturbation, so from now on we suppose such a change of coordinates to have been made.

Now eliminate the diagonal linear terms of $v$ by a shift of origin.  This is straightforward by completing the squares and produces a shift of order $\delta$.  Finally the effect of all the above on the constant terms is to make constant terms $(\tilde{\mu}_1, \tilde{\mu}_2)$ that are a near-identity diffeomorphism from $(\mu_1,\mu_2)$. 

Thus we are left with
\begin{eqnarray}
\dot{X}_1 &=& \tilde{\mu}_1 + X_1^2+\delta_1 X_2 + \eps_1 X_2^2 + HOT_1\label{eq:expansion} \\
\dot{X}_2 &=& \tilde{\mu}_2 + X_2^2 + \delta_2 X_1 + \eps_2 X_1^2 + HOT_2 \nonumber
\end{eqnarray}
for some coefficients $\delta_j$ and $\eps_j$ of order $\delta$, where $HOT_j$ denotes higher order terms of the form $HOT_j = f_j(X_j) + g_j(X)$ with $f_j = O(X_j^3)$ and $g_j = O(\delta |X|^3)$.

We make the assumption that $\delta_j \ne 0$, in fact of size similar to $\delta$ (though for some purposes, $\delta_j$ significantly larger than $\delta^{4/3}$ would suffice).  Then the $\eps_j$ terms are negligible relative to the $\delta_j$ terms.  The tidy way to deal with this is to push the $\eps_j$ terms into the $HOT_j$, so we shall consider this done.

The signs of $\delta_j$ will play a crucial role.  
We say the coupling coefficient $\delta_j$ is {\em excitatory} if $\delta_j>0$, {\em inhibitory} if $\delta_j<0$, by loose analogy with neuroscience. If both $\delta_j$ have the same sign we say the coupling is {\em mutualistic} ({\em mutually excitatory} if both are positive, {\em mutually inhibitory} if both are negative).  If the $\delta_j$ have opposite signs we say the coupling is {\em mixed}. 

Gap junction coupling gives the mutually excitatory case $\delta_1,\delta_2>0$, because its effect is to add current $I=\frac{V_1-V_2}{R}$ from neuron 1 to neuron 2, where $V_j$ are their electric potentials and $R$ is the resistance of the junction (taking a linear model).  This adds perturbation terms $\dot{V}_2=\frac{V_1-V_2}{RC_2}$, $\dot{V}_1=\frac{V_2-V_1}{RC_1}$, where $C_j$ is the capacitance of neuron $j$. Passage through the near resting state of a neuron corresponds to increasing $V$, so the angle coordinate $x_j$ for a SNIC neuron is oriented the same way as $V_j$ near the saddle-node.  
Up to a shift of origin of $V_j$, we have $x_j\sim\frac{V_j}{K_j}$ for a positive scale factor $K_j$ to make the quadratic coefficient equal 1.  So we read off that $\delta_1=\frac{K_2}{RC_1K_1}$, $\delta_2=\frac{K_1}{RC_2K_2}$, which are both positive.

By reversing time and $X_1$ if necessary, we can always take $\delta_1>0$, but one must remember to reverse time at the end of the analysis, which interchanges attractors and repellors for example (we will do this in Section~\ref{sec:MutInhibCase} for example). 

We shall switch notation back from $\tilde{\mu}_j$ and $X_j$ to $\mu_j$ and $x_j$, but it should be remembered that these are related to the original parameters and coordinates by a near-identity diffeomorphism.

As a simple example, consider the family
\begin{eqnarray*}
\dot{X}_1&=&\lambda_1-\cos X_1+\epsilon_1\sin X_2\\
\dot{X}_2&=&\lambda_2-\cos X_2+\epsilon_2\sin X_1
\end{eqnarray*}
on $\R^2/(2\pi\Z)^2$.  Then for coupling parameters $\epsilon_1,\epsilon_2=0$, the first equation has SNICs for $\lambda_1=\pm1$ at $X_1=0,\pi$ respectively, and the second equation has SNICs for $\lambda_2=\pm1$ at $X_2=0,\pi$ respectively.  Our coordinate and parameter changes for the resulting four cases are just shifts and scale changes.  Our special parameters $(\delta_1,\delta_2)$ for the four cases $\lambda_1=\pm1,\lambda_2=\pm1$ are just $(\lambda_1\epsilon_1,\lambda_2\epsilon_2)$.

\subsection{Reduced system}   
\label{sec:ReducedSyst}

In the triangle in parameter space (Figure \ref{fig:outside}) and the box $B$ ($|x_j|\le \eta$) in state space we study the approximate vector field $\hat{v}$
\begin{eqnarray}
\dot{x}_1 &=& \mu_1 + x_1^2 + \delta_1 x_2 \label{eq:approxvf}\\
\dot{x}_2 &=& \mu_2 + x_2^2 + \delta_2 x_1. \nonumber
\end{eqnarray}
Although the neglected higher order terms are small compared to $\mu$ and $x^2$, they are not necessarily small compared to $\delta_1 x_2$ and $\delta_2 x_1$, so one might ask why we retain the latter.  The idea is that all results inside the sub-box $|x| \ll \sqrt{\delta}$ will be accurate (because the higher order terms are dominated by the $\delta_1 x_2$ and $\delta_2 x_1$ terms there) and the results obtained outside this sub-box will turn out insensitive to the higher order terms anyway, because the dominant ones are of the form $\alpha_j x_j^3$ in $\dot{x}_j$ (rather than general cubics in both variables).

\begin{figure}[h] 
   \centering
   \includegraphics[width=2.8in]{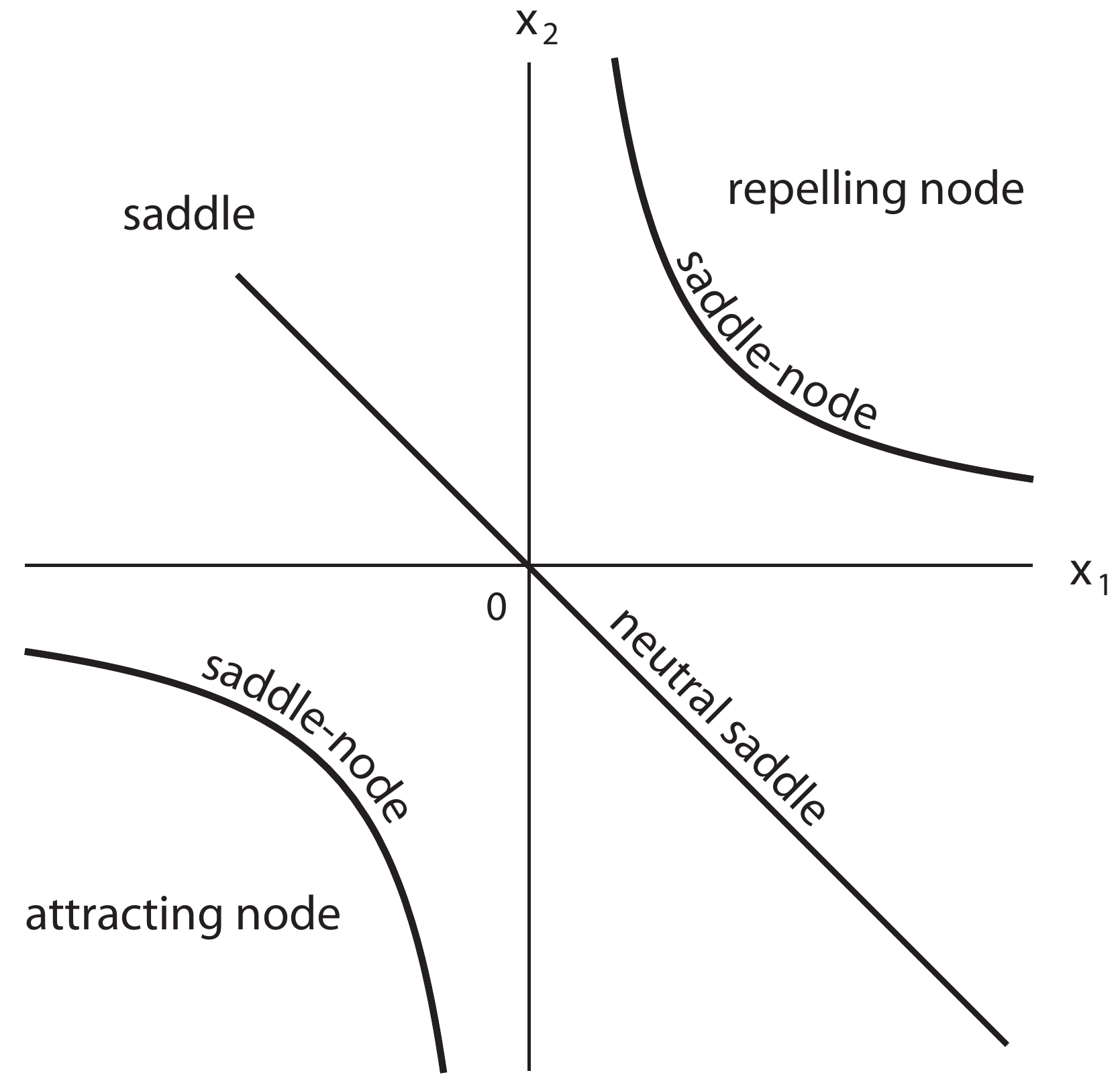} 
   \caption{Curves where $\det=0$ or $\tr=0$ on the manifold of equilibria, considered as a graph over $(x_1,x_2)$ for the mutualistic case $\delta_1\delta_2>0$.}
   \label{fig:equa}
\end{figure}

Its equilibria form the graph of a function from state space to parameter space:
\begin{eqnarray}
\mu_1 &=& -x_1^2 - \delta_1 x_2 \label{eq:equa} \\
\mu_2 &=& -x_2^2 - \delta_2 x_1. \nonumber
\end{eqnarray}
The type of the equilibrium is determined by the determinant and trace of the derivative $D\hat{v}$ of the vector field:
\begin{eqnarray}
\det D\hat{v} &=& 4x_1x_2-\delta_1\delta_2 \\
\tr \,D\hat{v} &=& 2(x_1+x_2) . \nonumber
\end{eqnarray}
In particular, the equilibrium is a saddle for $\det<0$, a sink for $\det>0,\tr<0$, and a source for $\det>0, \tr>0$.  The two branches of hyperbola in $(x_1,x_2)$ where $\det=0$ correspond to saddle-node equilibria, and the line where $\tr=0$ to neutral saddle or Hopf bifurcation (according as $\det <0, >0$).
Figure~\ref{fig:equa} illustrates the mutualistic case $\delta_1\delta_2>0$.  Note that $\tr=0$, $\det>0$ is impossible in this case so the $\tr=0$ curve is all neutral saddle.

The approximations for these curves are good for $|x| \ll \sqrt{\delta}$.  In fact they are good for all $|x| \le \eta\ll 1$. To see this, the dominant correction to the equation for a saddle-node equilibrium is
$$4x_1x_2-\delta_1\delta_2 = -6x_1x_2(\alpha_1x_1+\alpha_2x_2).$$
Because of the factor $x_1x_2$ on the right, the correction has relative size $O(\eta)$ and the saddle-node curves perturb to $4x_1x_2 = \delta_1\delta_2 (1+O(\eta))$.
The $\tr=0$ curve deforms to approximately $x_1+x_2 = -\frac{3}{2}(\alpha_1x_1^2+\alpha_2x_2^2)$, which is a shift of at most $O(\eta^2)$.

The projection of the manifold (\ref{eq:equa}) of equilibria to parameter space has fold curves where the determinant is zero.  To find their images in parameter space it is convenient to parametrise the two branches of hyperbola.  Here we make a separation of the analysis into the mutualistic and mixed cases.

In this paper we study the mutualistic case.  We will treat the mixed case in Paper II.

\section{The mutualistic case}\label{sec:analequ}    

The mutualistic case is $\delta_1\delta_2>0$.  By reversing time and the orientation of $x_1$ and $x_2$ if necessary, we can take $\delta_1$ and $\delta_2>0$.

\subsection{Analysis of equilibria}  

\begin{figure}[h] 
  \includegraphics[width=2in]{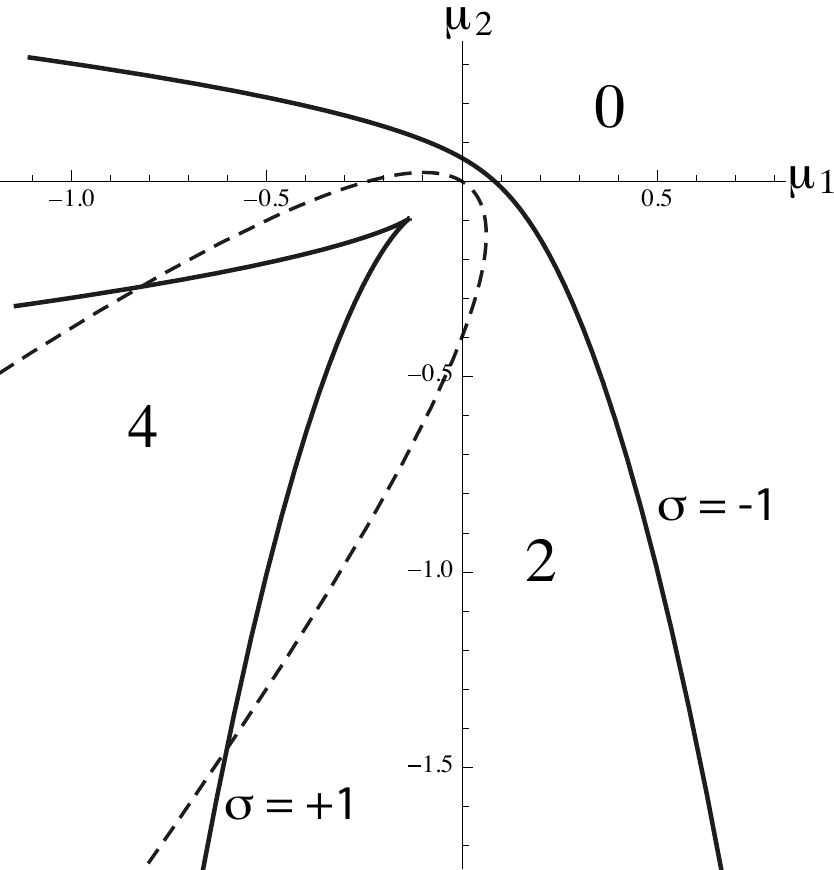} 
   \caption{Curves of saddle-node equilibrium (full) and of neutral saddle (dashed) in parameter space for $\delta_1\delta_2>0$ (drawn for $\delta_1=0.5$, $\delta_2=0.3$), also indicating the number of equilibria in the regions they separate.}
   \label{fig:sne}
\end{figure}

In the mutualistic case, we parametrise the saddle-node curves by 
\begin{eqnarray}
x_1 &=& \frac{\sigma}{2}\sqrt{\delta_1\delta_2} e^\theta \\
x_2 &=& \frac{\sigma}{2}\sqrt{\delta_1\delta_2}e^{-\theta}, \nonumber
\end{eqnarray}
with $\sigma = +1$ for the positive branch of the hyperbola of Figure \ref{fig:equa} and $-1$ for the negative branch.
It follows from \eqref{eq:equa} that the saddle-node curves project to
\begin{eqnarray}
\mu_1 &=& -\frac{\delta_1\delta_2}{4}e^{2\theta}-\frac{\sigma}{2}\delta_1\sqrt{\delta_1\delta_2}e^{-\theta} \\ \label{eq:sncurves}
\mu_2 &=& -\frac{\delta_1\delta_2}{4}e^{-2\theta}-\frac{\sigma}{2}\delta_2\sqrt{\delta_1\delta_2}e^\theta, \nonumber
\end{eqnarray}
which are drawn in Figure~\ref{fig:sne}.

Note that $d\mu_j/d\theta$, $j=1,2$, are non-zero along the $\sigma=-1$ curve, which we call the {\em outer saddle-node curve}, $\mu$ moving from lower right to upper left as $\theta$ increases.  They have a common zero along the $\sigma=+1$ curve, however, at $\theta=\theta_c$ where 
$$e^{3\theta_c} = \sqrt{\delta_1/\delta_2},$$
 which causes it to have a cusp at $$\mu_1=-\frac{3}{4}\delta_1^{4/3}\delta_2^{2/3},\ \mu_2 = -\frac{3}{4}\delta_1^{2/3}\delta_2^{4/3},$$ so we call it the {\em cusped saddle-node curve}.   The position of the corresponding degenerate equilibrium in state space then is $$x_1 = \frac12\delta_1^{2/3}\delta_2^{1/3},\ x_2 = \frac12\delta_1^{1/3}\delta_2^{2/3}.$$
In the case $\delta_1,\delta_2>0$ we are considering here, it is topologically a saddle but with zero exponent in the contracting direction.  It can be useful to rewrite the position of the cusp as
$$\mu_1 = -\frac34 \delta_1\delta_2 e^{2\theta_c},\ \mu_2 = -\frac34 \delta_1\delta_2 e^{-2\theta_c}, \ x_1 = \frac12 \sqrt{\delta_1\delta_2} e^{\theta_c},\ x_2 = \frac12 \sqrt{\delta_1\delta_2} e^{-\theta_c}.$$

The approximations for the saddle-node curves are good for all $|\mu| \le C\delta$, because the effect of the higher order terms is negligible anyway in $|\mu| \ll \delta$, and for the pieces of saddle-node curve in $-C\delta \le \mu_1 \le -K\delta^{2}$, $K$ large positive, where $\mu_1 \approx -x_1^2$, $\mu_2 \approx -\delta_2x_1$, we obtain 
$$\mu_2 = -\sigma  \delta_2\sqrt{-\mu_1}(1+O(\eta)).$$
The analogous result holds for the pieces in $-C\delta \le \mu_2 \le -K\delta^{2}$. 
In particular, the cusp in the saddle-node curve is stable to small perturbation of a two-parameter family, so survives.

Readers versed in singularity theory will recognise the saddle-node curves in Figure~\ref{fig:sne} as a slice through the unfolding of a hyperbolic umbilic singularity \cite{GG}.  Indeed, define a mapping $\varphi$ from $(x_1,x_2,\delta_1,\delta_2)$ to $(\mu_1,\mu_2,\delta_1,\delta_2)$ by \eqref{eq:equa}. Then its singularities (points where the rank of $D\varphi$ is less than 4) correspond to saddle-node equilibria  and the origin to a hyperbolic umbilic point. Because $\varphi$ is a stable mapping, all small smooth enough perturbations of the 4-parameter family of vector fields \eqref{eq:approxvf} have set of equilibria and singularity set smoothly equivalent  to that of  \eqref{eq:approxvf}.  Any smooth 2-parameter family of vector fields near the case $\delta_1=\delta_2=0$ has bifurcations of the set of equilibria  given by a slice through the unfolding of the hyperbolic umbilic.

Also shown in Figure~\ref{fig:sne} is the projection of the curve of neutral saddles to parameter space (Hopf bifurcation does not occur in the case $\delta_1\delta_2>0$), which is easily computed to be a parabola 
$$(\mu_1-\mu_2)^2+(\delta_1+\delta_2)(\delta_2\mu_1+\delta_1\mu_2) = 0.$$  
This is accurate, however, only for $|\mu| \ll \delta$; the effects of higher order terms can shift the curve by order $|\mu|^{3/2}$, which becomes the same size as the distance $(\delta_1+\delta_2)\sqrt{|\mu|}$ between the two sides of the parabola when $|\mu|$ approaches $C\delta$, so even allowing the possibility that the two sides cross.

It is useful to calculate the position of the two other equilibria for parameter values on the cusped saddle-node curve.  By factoring out the known double root, they are found to be at
\begin{eqnarray}
x_1 &=& -\frac{1}{2}\sqrt{\delta_1\delta_2}e^\theta+\eps\delta_1^{3/4}\delta_2^{1/4}e^{-\theta/2}\\
x_2 &=& -\frac{1}{2}\sqrt{\delta_1\delta_2}e^{-\theta} + \eps \delta_1^{1/4}\delta_2^{3/4}e^{\theta/2} \nonumber
\end{eqnarray}
for $\eps= \pm 1$, with $\eps=+1$ giving a saddle, $-1$ a sink.  

Furthermore, one can calculate the eigenvalues and eigenvectors for parameters on the cusped saddle-node curve (and on the outer saddle-node curve, so we give both cases).  The second eigenvalue of the saddle-node is $\tr = 2\sigma \sqrt{\delta_1\delta_2} \cosh\theta$ which is positive on the cusped curve, and the null eigenvector has slope $-\sigma\sqrt{\delta_2/\delta_1}e^\theta$ (by the {\em slope} of a vector $v = (v_1,v_2)$ we mean $v_2/v_1$).  A useful trick to save work is to note that the derivative $Dv$ is symmetric with respect to the inner product 
$$\langle \xi, \zeta \rangle = \xi_1\zeta_1/\delta_1 + \xi_2\zeta_2/\delta_2,$$ 
so the eigenvectors are perpendicular in this inner product.  This implies that the product of their slopes is $-\delta_2/\delta_1$.  So for example, the slope of the second eigenvector of the saddle-node is $\sigma\sqrt{\delta_2/\delta_1}e^{-\theta}$.

The eigenvalues $\lambda$ and slopes $s$ of the eigenvectors of any equilibrium are given by the following expressions:
\begin{equation}
\lambda = x_1+x_2 \pm\sqrt{(x_2-x_1)^2+\delta_1\delta_2}\label{eq:eval}
\end{equation}
\begin{equation}
s = \frac{1}{\delta_1} \left(x_2-x_1\pm\sqrt{(x_2-x_1)^2+\delta_1\delta_2}\right).\label{eq:slope}
\end{equation}
For the saddle which coexists along the cusped saddle-node curve, substitute the following expressions into \eqref{eq:eval} and \eqref{eq:slope}:
$$x_1+x_2 = \delta_1^{3/4}\delta_2^{1/4}e^{-\theta/2} + \delta_1^{1/4}\delta_2^{3/4}e^{\theta/2}  -\sqrt{\delta_1\delta_2}\frac{(e^\theta+e^{-\theta})}{2} = \sqrt{\delta_1\delta_2}\left(2 \cosh\frac{\theta-3\theta_c}{2} - \cosh\theta\right)$$
$$x_2-x_1 = \delta_1^{1/4}\delta_2^{3/4}e^{\theta/2} - \delta_1^{3/4}\delta_2^{1/4}e^{-\theta/2} +\sqrt{\delta_1\delta_2}\frac{(e^\theta-e^{-\theta})}{2} = \sqrt{\delta_1\delta_2} \left(2\sinh\frac{ \theta-3\theta_c}{2} + \sinh{\theta}\right).$$

\begin{figure}[h] 
   \centering                                            
   \includegraphics[width=3in]{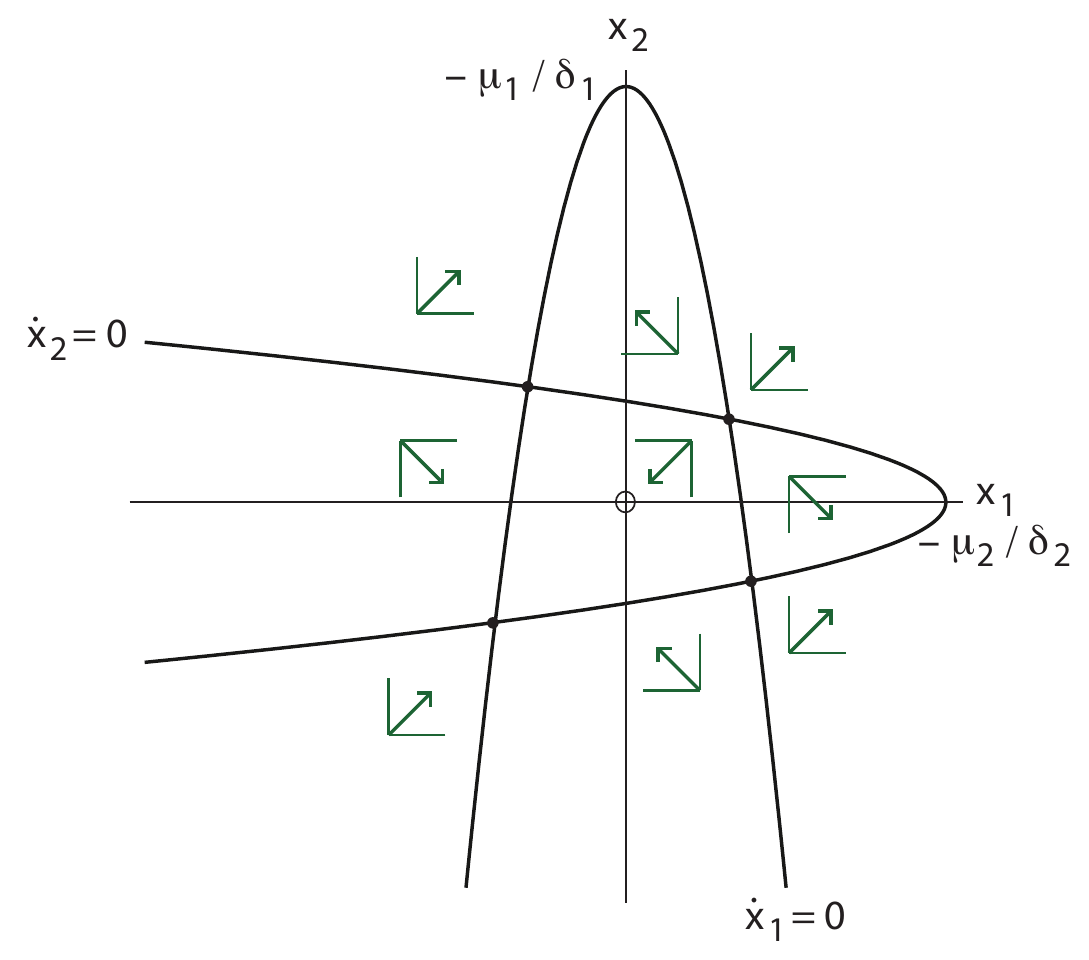} 
   \caption{Nullclines for parameters inside the cusped saddle-node curve ($\delta_1, \delta_2>0$).}
   \label{fig:nullclinesDelta2positif}
\end{figure}

\begin{figure}[h] %

   \centering                       
   \includegraphics[width=3in]{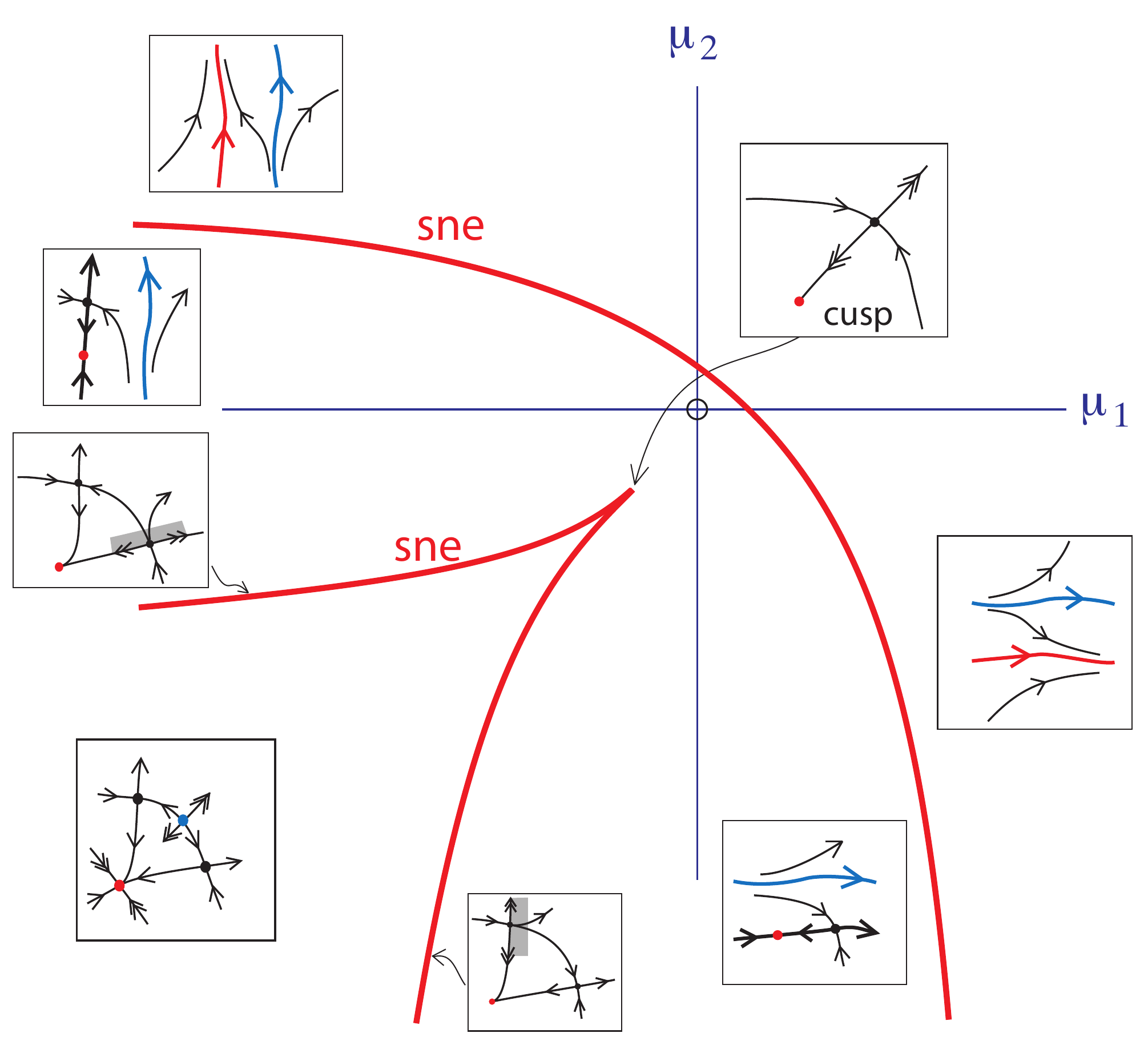} 
   \caption{Phase portraits in the box for the indicated parameter regimes ($\delta_1,\delta_2>0$). The winding ratio of the Poincar\'e and Cherry flows varies continuously between the extremes illustrated.  Shading near a saddle-node indicates its repelling half-plane.}
   \label{fig:PhasePortaitsBoxDelta2positif}
\end{figure}

The ways the equilibria connect within the box $B$ are deduced by studying the nullclines (e.g.~Figure~\ref{fig:nullclinesDelta2positif} for the region inside the cusped saddle-node of equilibria curve), and by computing the signs of the slopes of the eigenvectors at the equilibria from \eqref{eq:eval} and \eqref{eq:slope}.  In particular, for $\delta_1\delta_2>0$, the source and sink are always nodes (not foci) and for $\delta_1>0$, $\delta_2>0$ their fast and slow directions are as indicated on Figure~\ref{fig:PhasePortaitsBoxDelta2positif}.  Global connections will be analysed in subsections \ref{sec:intnbhd}
and \ref{sec:hetero}.

Thus the phase portraits in the box $B$ for the various parameter regimes, in particular as the parameters move along the cusped saddle-node curve, are as indicated on Figure~\ref{fig:PhasePortaitsBoxDelta2positif}. 

\subsection{Transit map}  
 \label{sec:transitmap} 

Next, to aid in understanding the global dynamics, we consider the transit map from $x_2=\eta$ to $x_2=L_2-\eta$, restricting attention to those trajectories that remain within the strip $|x_1|\le \eta$.  This analysis is independent of the sign of $\delta_2$, so logically could be done in the previous section, but it would have interrupted the presentation.

By the choice of $\eta\gg\sqrt{C\delta}$ and the restriction to the triangle in parameter space where in particular $|\mu_2|\le C\delta$, we have $\tilde{v}_2>0$ for $x_2 \in [\eta,L_2-\eta]$.  The transit map $x_1 \mapsto x'_1$ is given by integrating
\begin{equation*}
\frac{dx_1}{dx_2} = \frac{\tilde{v}_1}{\tilde{v}_2}
\end{equation*}
from $x_2 = \eta$ to $x_2 = L_2-\eta$.  As a first approximation, take $v_1= x_1^2$ in $|x_1|\le\eta$.
To this level of approximation the result is that 
\begin{equation}
\frac{1}{x_1}-\frac{1}{x'_1} = t_2 = \int_\eta^{L_2-\eta} \frac{dx_2}{\tilde{v}_2},\label{eq:x1'approx}
\end{equation}
where $t_2$ is the time taken, which may depend on $x_1$, but is dominated by the ends of the trajectory.  At the lower end $\tilde{v}_2 \approx x_2^2$ ($\mu_2$ is negligible here because $\eta \le x_2 \le L_2-\eta$), and at the upper end 
$\tilde{v}_2\approx(L_2-x_2)^2$ and so $t_2 \approx \frac{2}{\eta}$.  
As an explicit example, if 
\begin{equation}
v_2(x_2)=\left(\frac{L_2}{\pi}\right)^2 \sin^2\frac{\pi x_2}{L_2}\label{eq:v2}
\end{equation}
(which $\sim x_2^2$ for $x_2$ near $0$ and $\sim(L-x_2)^2$ for $x_2$ near $L_2$) then
$t_2=\frac{2\pi}{L_2}\cot\frac{\pi\eta}{L_2}\sim\frac{2}{\eta}$ for $\eta\ll L_2$.  

\begin{figure}[h] 
   \centering                                       
   \includegraphics[width=2.8in]{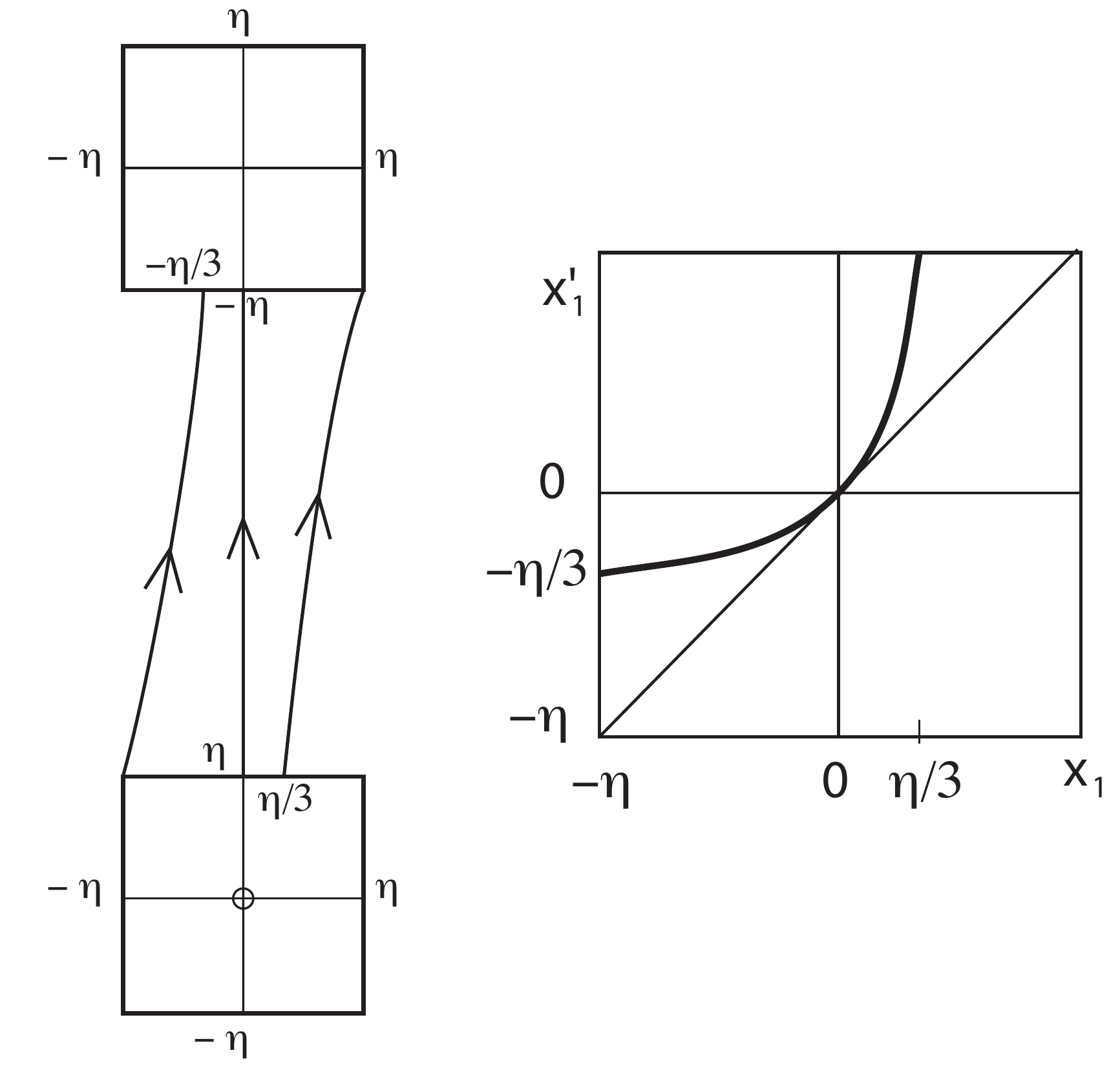} 
   \caption{Transit map}
   \label{fig:transitmap}
\end{figure}

Thus inserting $t_2\approx 2/\eta$ into \eqref{eq:x1'approx} the transit map is approximately
\begin{equation}\label{eq:x1'}
x'_1 = \frac{x_1}{1-2x_1/\eta} = \frac{\eta x_1}{\eta-2x_1}\,,
\end{equation}
valid for those transits remaining in the strip, i.e.~for $-\eta \le x_1 \le \eta/3$.  It takes the interval $[-\eta,\eta/3]$ to the interval $[-\eta/3,\eta]$ and has a degenerate fixed point at $0$.  The slope is bounded by $9$. See Figure~\ref{fig:transitmap}. 

The effects on the transit map of corrections to $v_1$ and $v_2$ at $\mu_1=\mu_2=0$ are shown in Appendix E to be  $O(x_1^2\log1/\eta)$ and the effects of parameters $\mu_j$ and the perturbation $\delta$ for $|\mu_j|\le C\delta$ are shown to be $O(\delta/\eta)$.

For the approximate transit map \eqref{eq:x1'}, 
$$x_1'-x_1=\frac{2x_1^2}{\eta-2x_1}.$$ 
It follows that the transit map moves all points to the right by at least $K\delta/\eta$ (some $K>0$) except when $|x_1| = O(\delta^{1/2})$. Within $|x_1| = O(\delta^{1/2})$ it moves points possibly left or right but by at most $K'\delta/\eta$, some $K'$.
The analogous result holds for the transit map from $x_1=\eta$ to $x_1=L_1-\eta$ in the strip $|x_2|\le \eta$.

\subsection{Extension of global dynamics into parts of the triangle}  
\label{sec:intnbhd}					

\begin{figure}[h] 
   \centering
   \includegraphics[width=4.3in]{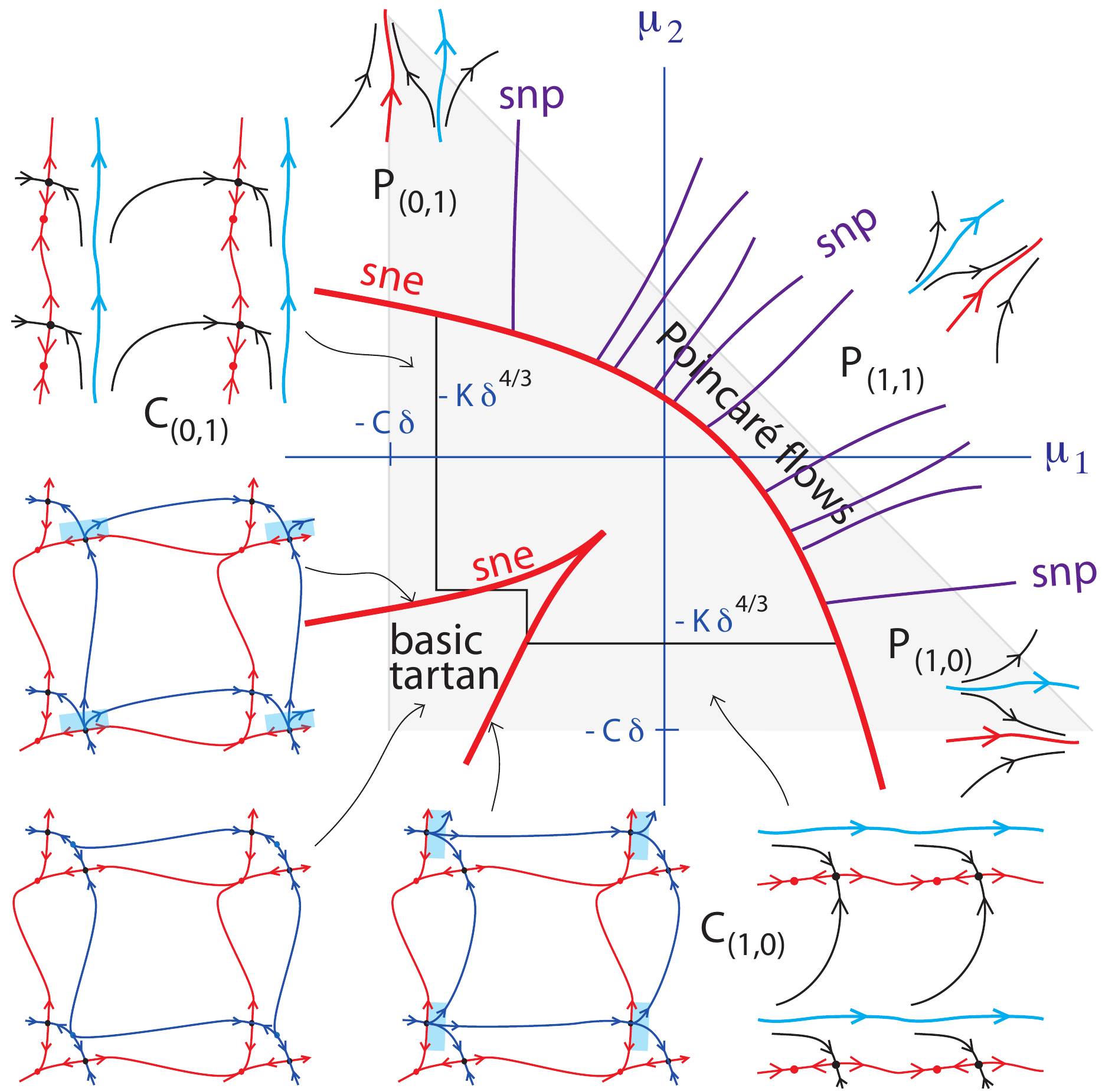} 
   \caption{Bifurcation diagram in parts of the triangle ($\delta_1,\delta_2>0$).}
   \label{fig:insidetriangle1}
\end{figure}

We can now describe the global dynamics in parts of the triangle in parameter space (Figure~\ref{fig:insidetriangle1}), leaving just a central region for analysis in subsection \ref{sec:hetero}.  

Firstly, there is $K$ large such that for $\mu_1\le-K\delta^{4/3}$ between the outer sne and lower cusped sne curves, the two vertical circles $C^\pm_1$ persist (though not necessarily $C^1$).  A sketch of the proof is given in Appendix F.

The invariant circles are periodic orbits above the outer saddle-node curve; the attracting one gains a saddle and sink on crossing the outer saddle-node of equilibrium curve and the repelling one gains a saddle and source on crossing the cusped saddle-node curve.  Below the cusped saddle-node curve, the study of the nullclines at the end of Subsection \ref{sec:ReducedSyst} shows that the equilibria connect in the box as in Figure~\ref{fig:PhasePortaitsBoxDelta2positif}.

Similarly, we obtain the analogous results in a strip along the horizontal boundary of the triangle.  In the intersection of the union of these boundary strips with the region inside the cusped saddle-node curve, the dynamics is fully mode-locked. 

Everywhere above the outer sne curve the dynamics consists of Poincar\'e flows, because there are no equilibria and there is a transverse section (e.g.~take $x_1=\frac{L_1}{2}$ where $\mu_1\ge \mu_2$ or $x_2=\frac{L_2}{2}$ where $\mu_2>\mu_1$).  The boundaries of the partial mode-locking strips are saddle-nodes of periodic orbits (snp).  In particular there is a curve of snp crossing the outer sne curve and delimiting the lower right zone of type $(1,0)$ Poincar\'e flow, and an analogous one for type $(0,1)$.  By the analysis at the beginning of this subsection, the type $(0,1)$ snp curve leaves the outer sne curve with $|\mu_1|\le K\delta^{4/3}$; but the most we can say for $\mu_2>-\mu_1$ is that it remains in $\mu_1\ge-C\delta$.  Similarly for the type $(1,0)$ snp curve.

It might be possible to obtain greater control over these snp curves as follows.  Consider type $(0,1)$, and take $\mu_2$ sufficiently positive.  The Lyapunov exponent of a periodic orbit in 2D is $\lambda = \oint \mbox{div}\, v\ dt$.  In the regime considered, the orbit is a graph $x_1(x_2)$, so this integral can be transformed to $\lambda = \int_0^{L_2} \frac{\text{div}\, v}{v_2}dx_2$.  The change in $x_1$ for an orbit segment making one revolution in $x_2$ is $\Delta x_1 = \int_0^{L_2} \frac{v_1}{v_2}\ dx_2$.  The conditions for an snp are $\lambda=0$, $\Delta x_1= 0$.  In the uncoupled case this happens at $x_1=0$, $\mu_1=0$.  The derivative of $\lambda$ with respect to initial condition $x_1$ on $x_2=0$ is $2T_2(\mu_2)$ which is large ($T_2(\mu_2) = \int_0^{L_2} \frac{dx_2}{v_2}$ is the period), so under perturbation $\lambda=0$ determines a nearby $x_1$.  The derivative of $\Delta x_1$ with respect to $\mu_1$ is $T_2$, so $\Delta x_1 = 0$ determines a nearby $\mu_1$.  Closer analysis, however, is required to determine more precise bounds on the snp curve.

Secondly, we are now able to give the phase portraits for the Cherry flows in between the outer and cusped sne curves near $\mu_1=-C\delta$ and $\mu_2=-C\delta$.  They are sketched in Figure~\ref{fig:insidetriangle1} for $\delta_1,\delta_2>0$, by analysis of the locations of the equilibria (cf. Figure~\ref{fig:PhasePortaitsBoxDelta2positif}).

Lastly, we prove that, writing 
$$\nu_j=-\mu_j\,,$$
the phase portrait is a basic tartan in the part of the region bounded by the cusped sne curve with both $\nu_j\gg \delta^2(\log{1/\delta})^4$.

\begin{figure}[htbp] 
   \centering					
   \includegraphics[height=3in]{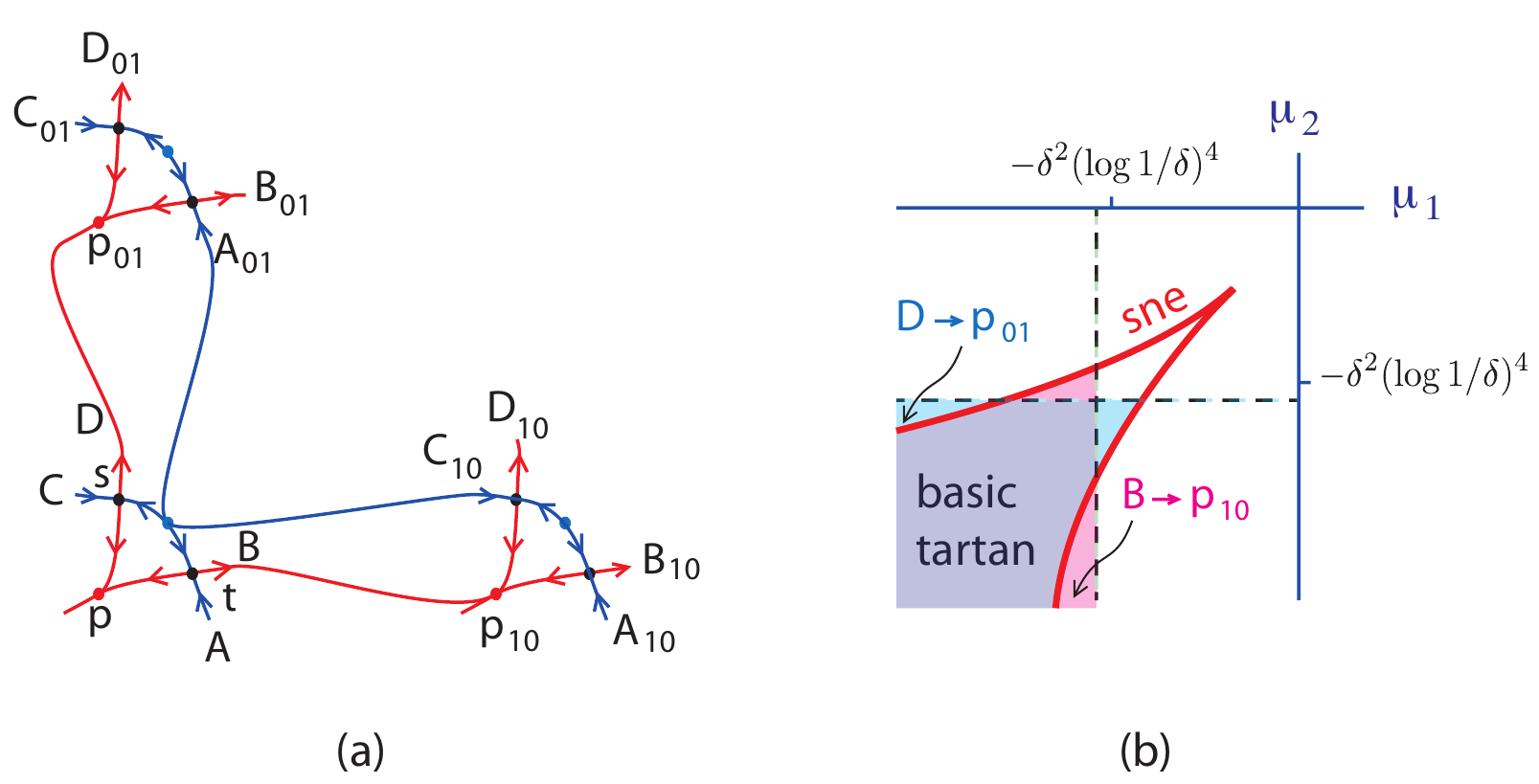}
   \caption{(a) Labelling of  the saddles $s$ and $t$, their relevant branches of local invariant manifolds $A,B,C,D$, and the sink $p$, also showing that the branch $D$ from saddle $s$ goes to the translated sink $p_{01}$ for $(\mu_1, \mu_2)$ not too close to the cusp, and similarly $B\to p_{10}$. (b) Regions of parameter space with the given connections.}
   \label{fig:ABCDlabels}  
\end{figure}

To begin, in the region bounded by the cusped saddle-node curve, we name the saddles $s$ and $t$ and the branches of the local invariant manifolds of the saddles $A, B, C, D$ as shown in Figure~\ref{fig:ABCDlabels}(a).  
On the part of the boundary going to the left from the cusp ($\theta > \theta_c$), $s$ becomes a saddle-node, with $A$ being in the null direction; on the part of the boundary going to the bottom ($\theta<\theta_c$), $t$ becomes a saddle-node, with $C$ being in the null direction.  At the cusp, $B$ merges with $D$.   Away from the cusp (where there is ambiguity), we may continue to use the labels for the branches of invariant manifold from the remaining saddle in the region between the two saddle-node curves (viz.~$C,D$ to the left, and $A,B$ to the bottom).  We denote the sink (which is in the negative quadrant) by $p$.  On the universal cover of the 2-torus we shall refer to the translation of $p$ by a vector $(m L_1,n L_2)$ with $m,n$ integer, as $p_{mn}$, and similarly for $s$, $t$, $A$, $B$, $C$, $D$.

For $\nu_2 \gg \delta^2 (\log{1/\delta})^4$, the branches $A$ and $D$ will be shown to be close to vertical for $|x_2|\le \eta = (\log{1/\delta})^{-1}$ (note the new choice of $\eta$), the transit map from $x_2=\eta$ to $L_2-\eta$ has a similarly small effect, but $D$ starts significantly to the left of $A$, so that $D$ reaches $x_2=L_2-\eta$ to the left of $A_{01}$. 

To prove the near-verticality of $A$ and $D$, note that they can be defined by
\begin{equation}
\frac{dx_1}{dx_2} = \frac{v_1(x_1,x_2)}{v_2(x_1,x_2)},
\end{equation}
starting from the appropriate equilibrium point which we denote in generality by $(x_1^e,x_2^e)$, and integrating in the appropriate direction of $x_2$ (increasing for $D$, decreasing for $A$).  Let us consider the case of $D$.  Then $x_1$ considered as a function of $x_2$ (which is true for $\eta$ not too large) is a fixed point of the map $T$ defined by
\begin{equation}
T[x_1](x_2) = x_1^e + \int_{x_2^e}^{x_2} \frac{v_1}{v_2}(x_1(x_2'),x_2')\ dx_2' .
\end{equation}
The map $T$ is a contraction in supremum norm on Lipschitz functions from $[x_2^e,\eta]$ to $\R$ satisfying $x_1(x_2^e)=x_1^e$ with Lipschitz constant $K$, say, if $\eta$ is sufficiently less than $L_2/2$, so achieved for $\delta$ small enough.  Call its contraction constant $\lambda < 1$; it can be made as small as we want by suitable choice of $K$ and $\delta$.
To bound the distance of the fixed point of $T$ from the constant function $x_1(x_2)=x_1^e$ it suffices to estimate the distance of the image of the constant function from the constant function and divide it by $(1-\lambda)$.

Now $v_1(x_1^e,x_2)$ and $v_2(x_1^e,x_2)$ have a common factor $(x_2-x_2^e)$, and for $x_2$ small the ratio is approximately $\delta_1/(x_2+x_2^e)$.  This approximation is not particularly accurate for larger $x_2$ but is good enough.  Then
\begin{equation}
\int_{x_2^e}^{x_2} \frac{v_1}{v_2}(x_1^e,x_2')\ dx_2' \approx \delta_1 \log\frac{x_2+x_2^e}{2x_2^e} .
\end{equation}
It follows that 
\begin{equation}
x_1(\eta)-x_1^e \le k\delta_1\log\frac{\eta+x_2^e}{2x_2^e},\label{eq:x1-x1e}
\end{equation}
with $k$ a little larger than $(1-\lambda)^{-1}$.

A similar bound is obtained for $A$, with $x_2^e$ replaced by $|x_2^e|$.

The effect of the transit of $D$ from $x_2=\eta$ to $L_2-\eta$ is of order $\delta/\eta$ as in subsection \ref{sec:transitmap}.

Now to complete the analysis, the horizontal distance $\Delta x_1$ between the two equilibria for given $\nu_2$ is at least that for the worst case of $\nu_1$, namely on the lower branch of the cusped sne curve, because as $\nu_1$ increases from this, the $\dot{x}_1$ nullcline rises, thereby separating the two equilibria further.  For this worst case, and taking $\nu_2$ significantly larger than $\delta^2$, the formulae for the equilibria give
\begin{equation}
\Delta x_1 \sim \delta_1^{3/4}\delta_2^{1/4} e^{-\theta/2},
\end{equation}
with
\begin{equation}
\nu_2 \sim \frac14 \delta_1 \delta_2 e^{-2\theta}.
\end{equation}
Thus
\begin{equation}
\Delta x_1 \sim \sqrt{2\delta_1} \nu_2^{1/4}.
\end{equation}

On the sne curve, we also have the approximation
\begin{equation}
|x_2^e| \sim \frac12 \sqrt{\delta_1 \delta_2} e^{-\theta} \sim \sqrt{\nu_2}
\end{equation}
for both equilibria (and this does not change much if $\mu_1$ is moved to the left).  So the shifts in $D$ and $A$ on reaching $x_2 = \pm\eta$ respectively are from \eqref{eq:x1-x1e} at most approximately $k\delta_1 \log\frac{\eta}{2\sqrt{\nu_2}}$.

Our goal is attained if
\begin{equation}
\delta_1 \log\frac{\eta}{ 2\sqrt{\nu_2}} + \frac{\delta}{\eta} \ll \sqrt{2\delta_1} \nu_2^{1/4}.
\end{equation}
Choosing $\eta = (\log{1/\delta})^{-1}$ (for which the value of the left hand side is of the same order as its minimum, yet $\eta \to 0$ as $\delta \to 0$), we obtain the sufficient condition
\begin{equation}
\nu_2 \gg \delta^2 (\log{1/\delta})^4,
\end{equation}
for $D$ to pass to the left of $A_{01}$. 
 
Similarly, $B$ passes under $C_{10}$ if $\nu_1 \gg \delta^2 (\log{1/\delta})^4$. Hence we obtain the global connections $D\to p_{01}$ and $B\to p_{10}$ in the regions indicated in Fig.~\ref{fig:ABCDlabels}(b).

\subsection{Heteroclinic connections and consequences}
\label{sec:hetero}

\begin{figure}[htbp] 
   \centering				
   \includegraphics[width=2.5in]{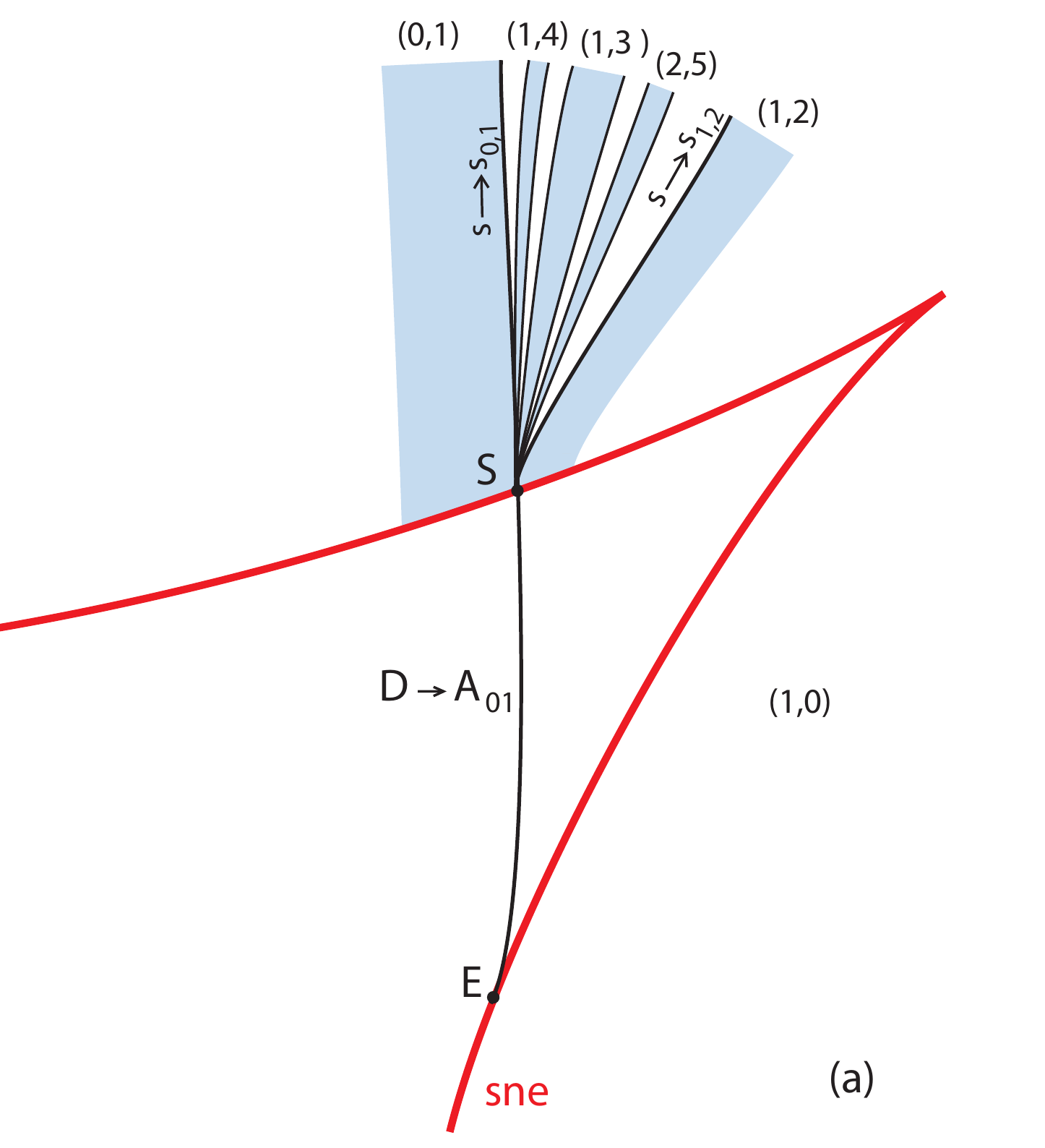} \hspace{1cm}   \includegraphics[width=2in]{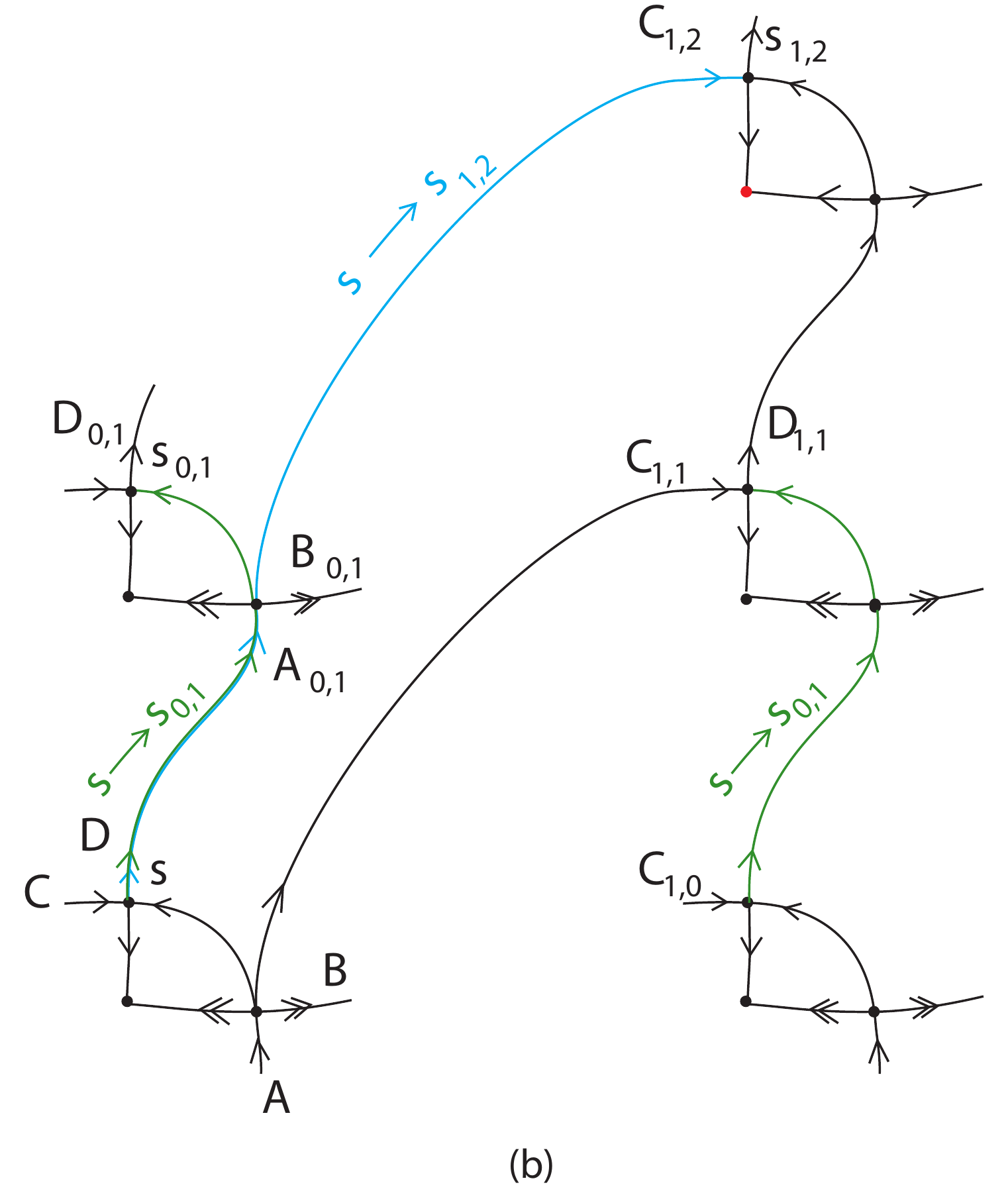} 
   \caption{(a) Curve of heteroclinic connection $D\to A_{01}$, $S$-point and resulting saddle-node fan between curves of rotational homoclinic connections $s\to s_{0,1}$ and $s\to s_{1,2}$ (in case $n=1$).  (b) Phase portrait at the $S$-point.}
   \label{fig:hetero1snfan}
\end{figure}

To determine the structure of the bifurcation diagram in the remaining central region of Fig.~\ref{fig:insidetriangle1}, we 
consider first the special case of systems which are symmetric with respect to simultaneous interchange of $x_1$ with $x_2$ and $\mu_1$ with $\mu_2$, thus in particular $\delta_1=\delta_2$ and $L_1=L_2$.  Then by symmetry at the cusp we have Cherry flow of type $(1,1)$, but for $\nu_2\gg\delta^2\left(\log\frac{1}{\delta}\right)^4$, $D$ passes to the left of $A_{01}$ as we just showed, so a curve of  heteroclinic connection $D \to A_{01}$ must occur, as indicated in Fig.~\ref{fig:hetero1snfan}(a), separating the cusp from the region $\nu_2\gg\delta^2\left(\log\frac{1}{\delta}\right)^4$.  The curve of heteroclinic bifurcation leaves the lower branch of cusped sne curve tangentially; this is a bifurcation not seen in \cite{BGKM,BM} so we give it a new name, $E$ point, and analyse its generic unfolding in Appendix G.   Where the curve of heteroclinic bifurcation hits the leftward branch of the cusped sne curve, we obtain the phase portrait of Figure~ \ref{fig:hetero1snfan}(b).  This is
 called an $S$ point in  the terminology of \cite{BGKM}, and it produces a ``saddle-node fan'' of Cherry flows.  The range of winding ratios produced in the saddle-node fan depends on the copy of the sink $p$ to which $B$ connects there.  If $B \to p_{1n}$ then the saddle-node fan generates all winding ratios from $(0,1)$ to $(1,n+1)$.  See Figure \ref{fig:hetero1snfan}(a) for the case $n=1$. 

By reflection symmetry there is an analogous curve of heteroclinic connection $B \to A_{10}$, generating a saddle-node fan of all winding ratios from $(1,0)$ to $(m+1,1)$ if $D\to p_{m1}$.   If there is no other bifurcation along the lower cusped sne curve between the $D \to A_{01}$ heteroclinic and this saddle-node fan, then $n=1$, and we shall suppose this case for the next five paragraphs.

\begin{figure}[htbp] 
   \centering				
   \includegraphics[width=5in]{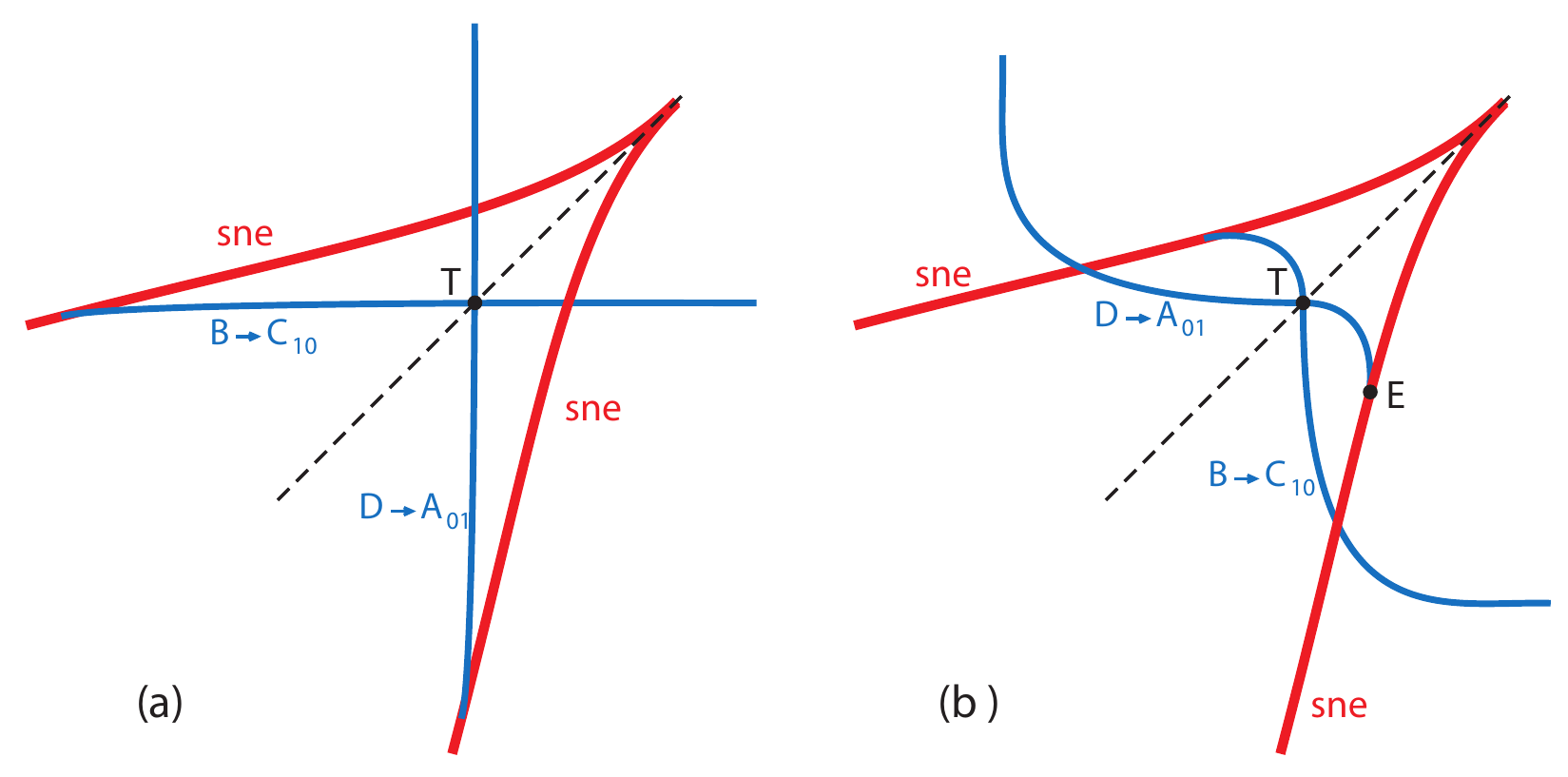}  
    \caption{Crossing  of curves of heteroclinic bifurcation makes a $T$-point on the diagonal in the symmetric case. Two possibilities: (a), (b).}
   \label{fig:TptSym}
\end{figure}

\begin{figure}[htbp] 
   \centering			
   \includegraphics[height=2.5in]{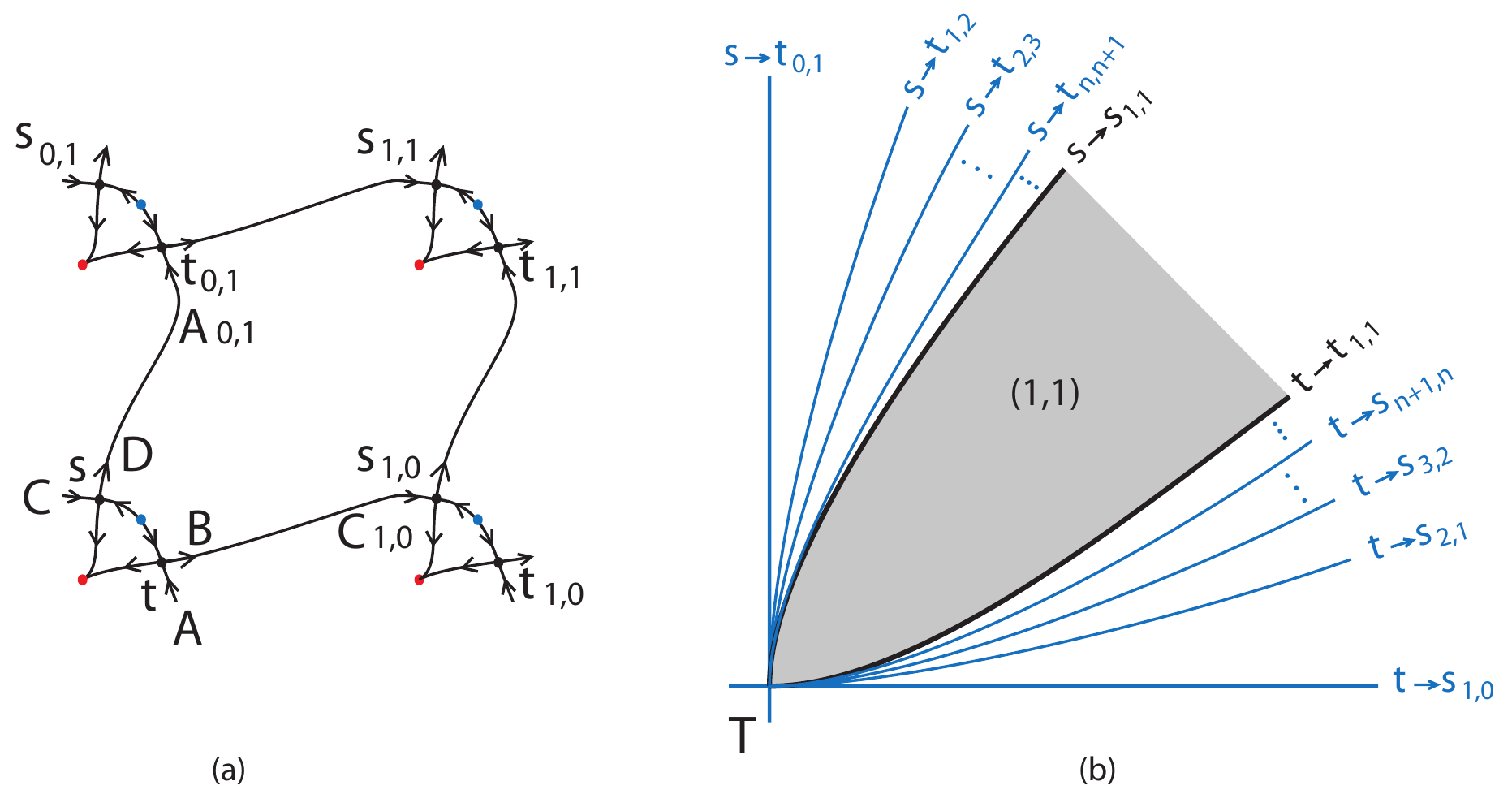} 
   \caption{$T$-point. (a) Phase Portrait. (b) Unfolding. The lines of heteroclinic connections that emanate from the $T$ point, $B\to C_{n+1,n}$  and $D\to A_{n,n+1}$, $n=1,2,...$, accumulate on lines of homoclinic connections $B\to A_{11}$, $D\to C_{11}$, respectively.}
   \label{fig:Tpoint}
\end{figure}

The curves of heteroclinic bifurcation $D \to A_{01}$ and $B \to C_{10}$ cross on the diagonal at a $T$ point (in the terminology of \cite{BGKM}).  Two ways this might happen are shown in Fig.~\ref{fig:TptSym} (others with more crossings can be envisaged).  We suspect only case (a) occurs but despite much effort to prove that the curves of heteroclinic connection are roughly vertical and horizontal we did not succeed.  Nevertheless, let us concentrate attention on case (a), and we will argue at the end that case (b) can not occur.  

The phase portrait is shown in Figure~\ref{fig:Tpoint}(a), and its unfolding produces a fan of curves of heteroclinic and rotational homoclinic bifurcation and a partial mode-locking tongue, as shown in Fig.~\ref{fig:Tpoint}(b). Note that in this symmetric case, $x_1+x_2=\delta_1$ for the two saddles, so both are repelling; this makes the exponents of the saddles less than 1 (for a saddle with eigenvalues $-\lambda < 0 < \mu$, the {\em exponent} $\alpha = \lambda/\mu$), so the bifurcation diagram is the time-reverse of the first case of Fig. ~4.23(c) of \cite{BGKM}.

\begin{figure}[htbp] 
   \centering			
   \includegraphics[width=2in]{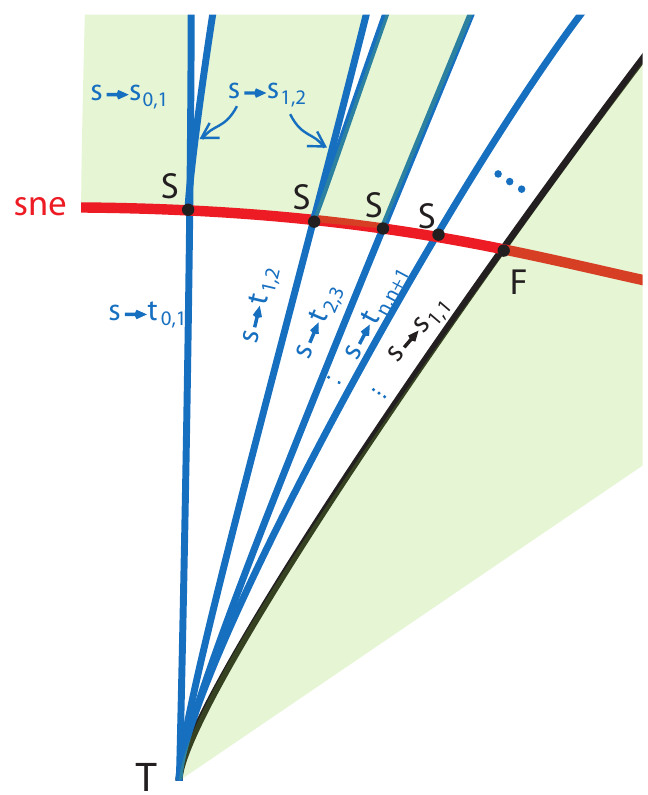} 
   \caption{A sequence of $S$ points accumulating onto an $F$ point.}
   \label{fig:Fpoint}
\end{figure}

\begin{figure}[htbp] 
   \centering			
   \includegraphics[width=4in]{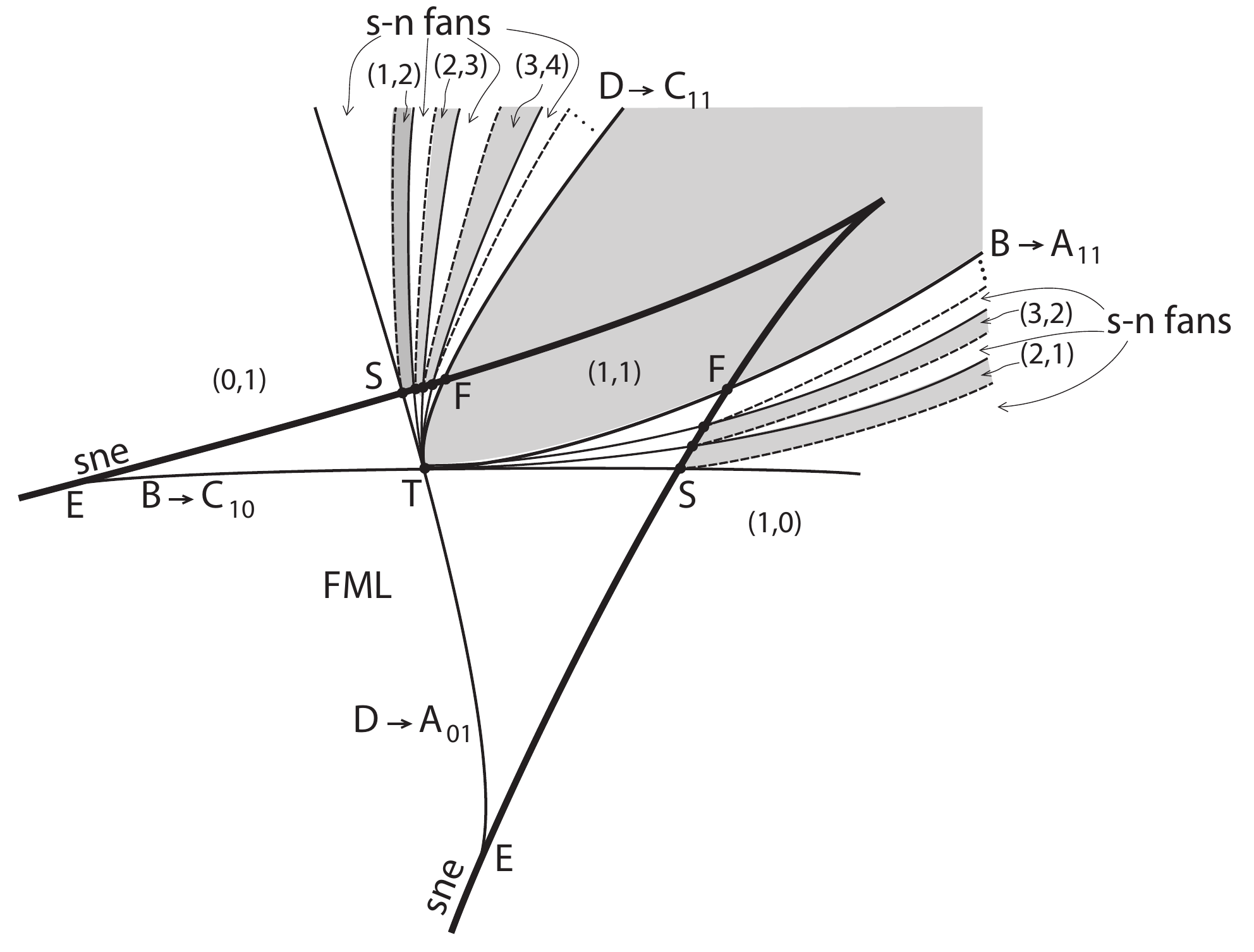}   
   \caption{Simplest bifurcation diagram around the cusp for $\delta_1,\delta_2>0$.  The lines of heteroclinic connections that emanate from the $T$ point cross the sne curve at saddle-node fan points ($S$ points) and continue as lines of homoclinic bifurcations, forming boundaries of regions of Cherry flow of type indicated. The dashed lines correspond to homoclinic bifurcations that use the other branches of the manifolds of the remaining saddle. }
   \label{fig:simplestdiag}
\end{figure}

Now if we suppose the bifurcation diagram near the $T$ point extends to the cusped sne curve, then we obtain a sequence of saddle-node fans from each of the indicated curves of heteroclinic connection (of which the first is the one we already found).  They accumulate at what we christen an $F$ point as shown in Fig.~\ref{fig:Fpoint}; it is an analogue of the $M$ point of \cite{BGKM} but with its saddle-node of periodic orbits replaced by a rotational homoclinic connection.  The resulting bifurcation diagram is sketched in Figure~\ref{fig:simplestdiag}.

\begin{figure}[htbp] 
   \centering            
   \includegraphics[height=3in]{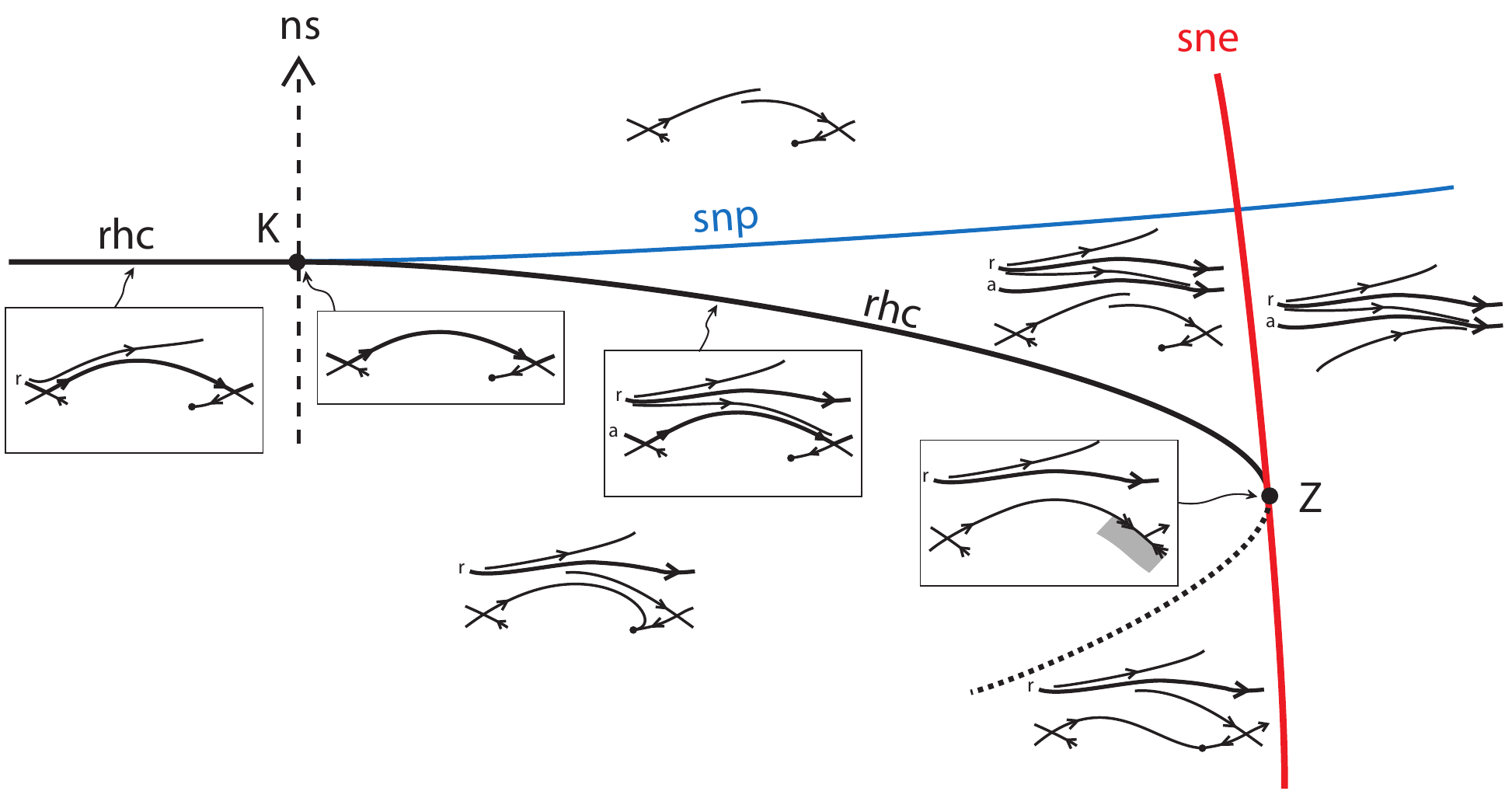} 
   \caption{A horn of coexistence of attracting equilibrium and periodic orbit, associated with $K$ and $Z$ points (cf. \cite{BGKM}). }
   \label{fig:KZhorn}
\end{figure}

Extending further in parameter space, where the resulting partial mode-locked tongues cross the curve of neutral saddle ($K$ points of \cite{BGKM}), their boundaries are replaced by curves of saddle-node periodic orbits, with a curve of rotational homoclinic connection between the $K$ point and a $Z$ point (saddle-node loop) on the outer sne curve.   We recall from \cite{BGKM}  in Figure~\ref{fig:KZhorn}, how $K$ and $Z$ points can produce a horn of coexistence of attracting periodic orbit with attracting equilibrium.  Calculation of the generic unfolding of $K$ points was given by Dumortier, Roussarie and Sotomayor \cite{DRS} and of $Z$ points by Schechter \cite{Sch}

\begin{figure}[htbp] 
   \centering             
   \includegraphics[width=4in]{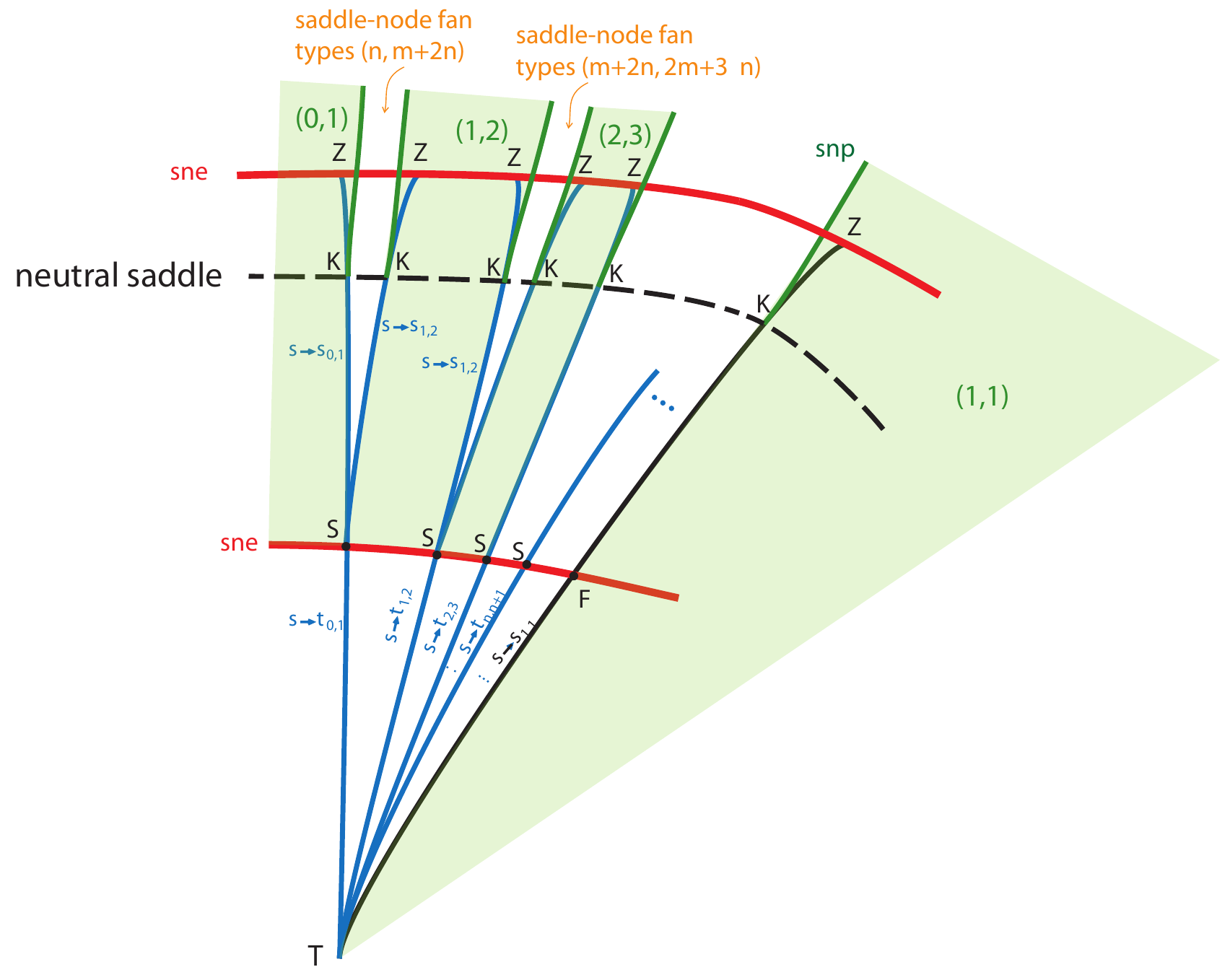} 
   \caption{Extending further in parameter space for $\delta_1,\delta_2>0$, where the resulting partial mode-locked tongues cross the curve of neutral saddle ($K$ points), their boundaries are replaced by curves of saddle-node periodic orbits, with a curve of rotational homoclinic connection between the $K$ point and a $Z$ point (saddle-node loop) on the outer sne curve. }
   \label{fig:creationofKpoints}
\end{figure}

The resulting bifurcation diagram is sketched in Fig.~\ref{fig:creationofKpoints}.

It might be that some of the curves of heteroclinic connection cross the curve of neutral saddle before reaching the cusped sne curve, which would make a slightly different picture.

For families that are not symmetric, but close to symmetric, the above conclusions continue to hold.  Simply the $T$ point is not necessarily on the diagonal.  Further away from symmetry, however, various changes can occur. In particular, the $T$ point could cross the neutral saddle curve of Figure~\ref{fig:sne}, changing the exponent of one saddle relative to 1, producing other cases of unfolding of the $T$ point (see \cite{BGKM}).  Or it could happen that the first bifurcation on sliding up the lower cusped sne is the saddle-node fan generated by $B \to C_{10}$ bifurcation, with a $D \to A_{01}$ connection happening higher up.  This would make the first saddle-node fan cover all winding ratios from $(1,0)$ to $(1,1)$ and necessitate other changes, but eventually there would be a $T$ point (with different homotopy type of heteroclinic cycle), generating a picture similar to the above (or one of its variants depending on the exponents of the saddles).

To complete this section, we conjecture that case (b) of Figure \ref{fig:TptSym} can not occur.  Our reason is that the $E$ point at the end of the $D\rightarrow A_{01}$ curve produces flow of type $(1,1)$ just outside the cusped region, whereas the $(1,1)$ tongue produced by the $T$ point is contained in the quadrant containing the cusp point, and it seems unlikely to us to see non-monotonicity of the winding ratio on turning round the cusp.

\section{Consequences for attractors} 

For applications the most important feature is attractors.  We can deduce what bifurcations happen to attractors in the mutually excitatory case by studying the bifurcation diagrams and phase portraits of Section~\ref{sec:analequ}. As mentioned, gap junction coupling gives the mutually excitatory case.  To deduce those for the mutually inhibitory case, we must time-reverse the flows and the orientations of $x_1$ and $x_2$.

\subsection{Mutually excitatory case} 

\begin{figure}[htbp] 
   \centering					
   \includegraphics[height=3in]{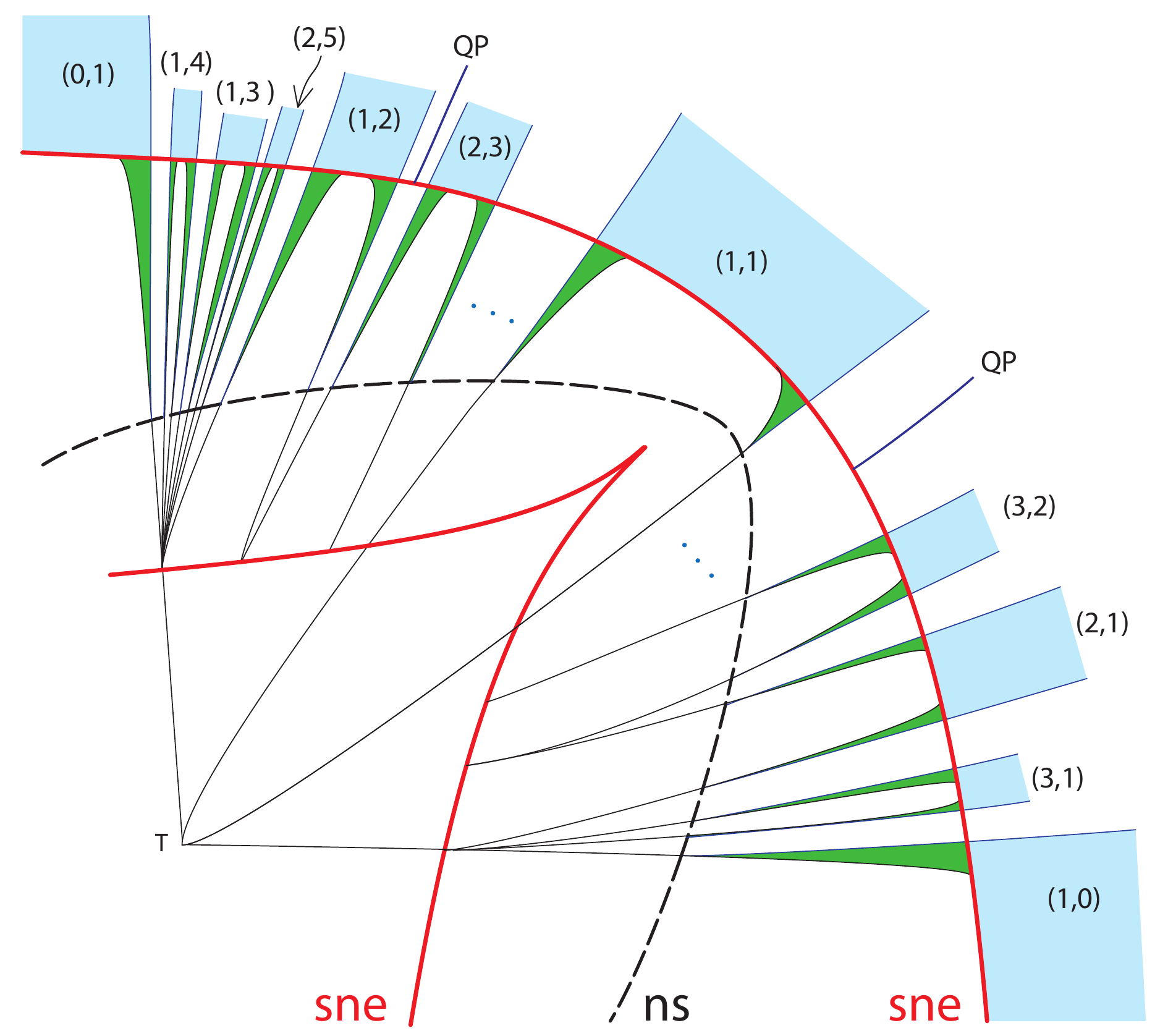}
   \caption{Attractors in the mutually excitatory case. There is an attracting equilibrium everywhere inside the outer saddle-node curve. In the indicated strips there is an attracting periodic orbit of the indicated homotopy type.  The strips extend across the outer sne curve into two horns where attracting equilibrium and periodic orbit coexist. Outside the outer sne curve the strips of periodic attractors are interleaved by curves on which the whole torus has quasi periodic flow of given homotopy type.}
   \label{fig:AttrMutExcitCase}  
\end{figure}

Corresponding to the bifurcation diagram of Figure~\ref{fig:creationofKpoints} for $\delta_1,\delta_2>0$, we obtain the types of attractors indicated in Figure~\ref{fig:AttrMutExcitCase}.  

In the whole of the region to the South-West of the outer sne curve there is an attracting equilibrium.  The region outside the outer sne curve is divided into strips where there are one (or more) attracting homotopically non-trivial periodic orbits.  This implies a periodic spiking pattern of the pair of oscillators.  There is one strip for each ratio of firing frequencies between 0 and $\infty$.  In between the strips are curves on which the whole torus  is an attracting quasiperiodic flow.  Each rational strip extends across the outer sne curve in two horns terminating on the neutral saddle curve (except just one horn for the $(0,1)$ and $(1,0)$ strips).  Thus in these horns two attractors coexist: an equilibrium and a homotopically non-trivial periodic orbit.  This produces hysteresis effects.  Deformations of Figure~\ref{fig:creationofKpoints} as mentioned towards the end of Subsection~\ref{sec:hetero} will produce deformation of Figure~\ref{fig:AttrMutExcitCase}  but no qualitative change in this description of the attractors.

\subsection{Mutually inhibitory case} \label{sec:MutInhibCase} 

\begin{figure}[htbp] 
   \centering					
   \includegraphics[height=3in]{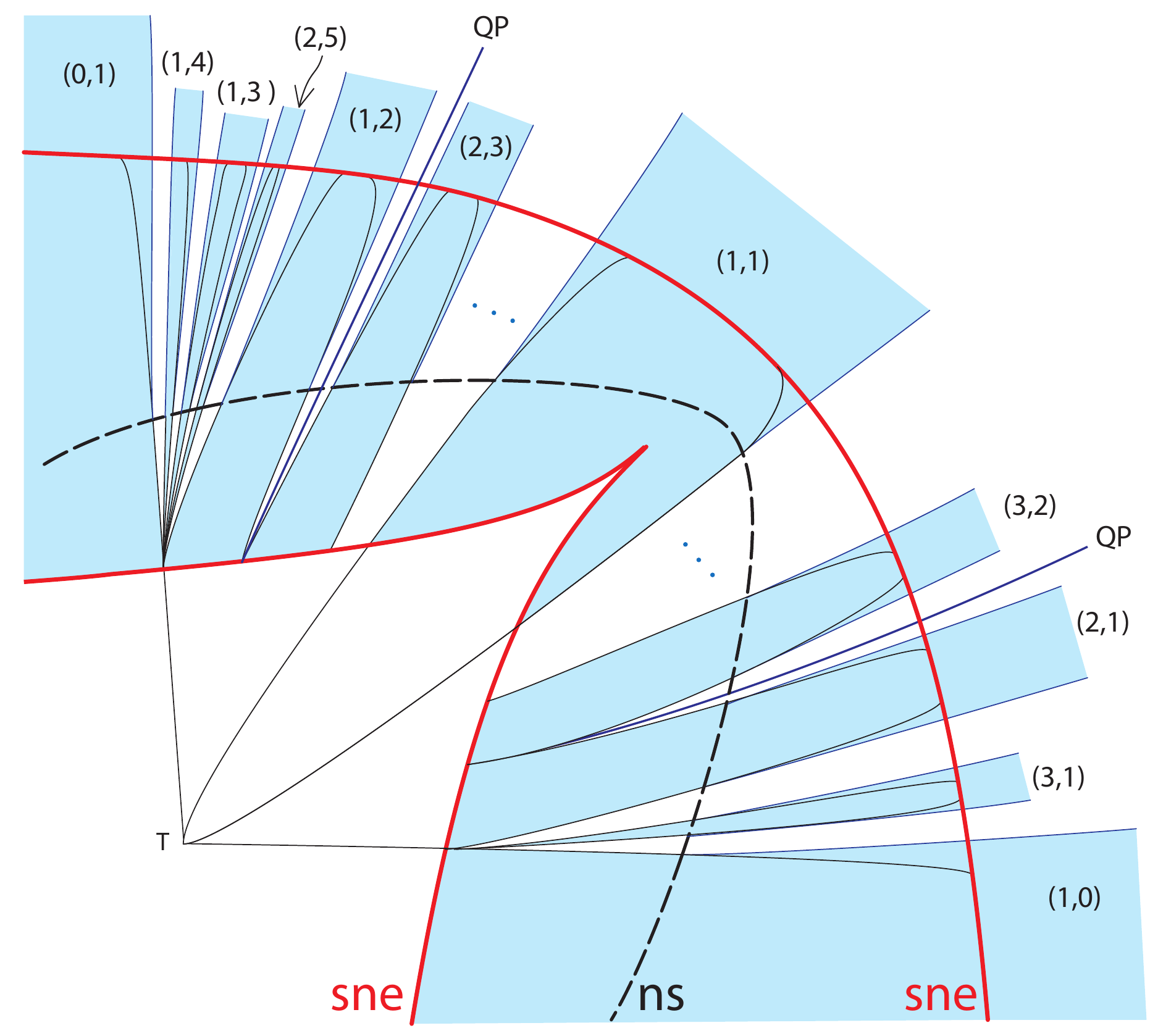}
   \caption{Attractors in the mutually inhibitory case. Inside the cusped saddle-node curve the only attractor is an equilibrium.  The region outside the cusped sne curve is divided into strips where there is an attracting homotopically non-trivial periodic orbit, interleaved by curves with a quasi-periodic attractor.}
   \label{fig:AttrMutInhibCase}  
\end{figure}

Time-reversing the results of Section~\ref{sec:analequ}, we obtain Figure~\ref{fig:AttrMutInhibCase} for the attractors corresponding to the bifurcation diagram of Figure~\ref{fig:creationofKpoints} in the case  $\delta_1,\delta_2<0$.
The situation is simpler  than the mutually excitatory one.  Inside the cusped sne curve there is an attracting equilibrium  and it attracts almost everything.  The region outside the cusped sne curve is divided into strips  with an attracting periodic orbit of type $(0,1)$, 
$(1,2)$, $(2,3)$,..., $(1,0)$, $(2,1)$, $(3,2)$,... and $(1,1)$ (which has the cusp cut out of it) and tongues for all other rationals, separated by curves of quasiperiodic  attractor. Outside the outer sne curve, the quasiperiodic attractor is the whole torus, but in between the two sne curves it is a ``cantorus'' (also known as Denjoy counterexample or Cherry attractor).

There is no region of coexistence of attractors, except that there may be more than one periodic orbit of given type.

The diagram may change significantly if Figure~\ref{fig:creationofKpoints} deforms in some of the ways described towards the end of Subsection~\ref{sec:hetero}.  Which rational gets the cusp and which are strips rather than tongues depends at least on the type of $T$ point.  But the overall result is roughly the same.

\section{Conclusion}  

We have analysed the generic interaction of two dynamical systems with saddle-node on an invariant circle (SNIC) bifurcation.   We identified two key coupling parameters $\delta_1$ and $\delta_2$.  We studied the ``mutualistic'' case $\delta_1\delta_2>0$ in detail.  We found that even the simplest sub-cases  of this have fairly complicated bifurcation diagrams.   We found that the attractors can be equilibrium, homotopically non-trivial periodic orbit, quasi-periodic torus or quasiperiodic cantorus and that in the mutually excitatory case, which is expected to be the relevant case for gap junction coupling, there are parameter ``horns'' where attracting equilibrium and periodic orbit coexist.  The results are expected to be a useful guide to the study of the dynamics of networks of class I neurons, Josephson junctions  and other situations where SNIC bifurcations occur.

Results for the mixed case $\delta_1\delta_2<0$, which turns out to be even more complicated, will be presented in a separate paper.  

\section*{Acknowledgements}
We are grateful to Steve Coombes for confirming that this analysis would be of interest to neuroscientists and to the University of Warwick for study leave and the Universit\'e Libre de Bruxelles for hospitality which enabled us to write the bulk of the paper.  We thank the referees for their constructive suggestions for improvement.  The results were presented at a minisymposium in Equadiff 2011 supported by the Engineering and Physical Science Research Council (grant number EP/G021163/1).

\section*{Appendix A: Coordinate change to make $b(\mu)=1$}
We give here an example of a parameter-dependent coordinate change $X(x,\mu)$ on a circle, preserving the length $L$, which makes the coefficient $b(\mu) = 1$ in (\ref{eq:sne}) for all small $\mu$.  

Let $$X(x,\mu) = x - \alpha(\mu)\sin{kx},$$ with $k = \frac{2\pi}{L}$ and $|\alpha| < \frac{1}{k}$ for invertibility.  Then $$\dot{X} = (1-k\alpha\cos{kx})(a(\mu) + b(\mu) x^2 + o(x^2))$$ and $$X = (1-k\alpha)x + O(x^3),$$ so $$x=X/(1-k\alpha) + O(X^3)$$ and $$\dot{X} = (1-k\alpha(1-\frac{k^2x^2}{2}))(a + bx^2) + o(x^2) = (1-k\alpha)a + \frac{\frac12{k^3\alpha a} + (1-k\alpha) b} {(1-k\alpha)^2} X^2 + o(X^2).$$  Now we choose $\alpha(\mu)$ to make the coefficient of $X^2$ equal $1$.  This is a quadratic equation for $\alpha$ with a simple root near $0$ because $b$ is near $1$.  

Thus we obtain $\dot{X} = A(\mu) + X^2 + o(X^2)$ for some function $A$, with $A(0)=0$.

\section*{Appendix B: Coordinate change to make $\frac{\partial v}{\partial \mu}\ge c$}

We prove here that for a SNIC $\dot{x} = v(x,\mu)$ with form \eqref{eq:dotxmux2} and for $c\in(0,1)$, by a parameter-dependent coordinate change $X(x,\mu)$ preserving the length $L$ of the circle, we can make $\frac{\partial v}{\partial \mu} \ge c$ for all $x$ for $\mu$ small.

The vector field in the new coordinate $X$ is given by $\dot{X} = X_x(x,\mu) v(x,\mu)$ (subscript denoting partial derivative), with $x(X,\mu)$ determined by $X = X(x,\mu)$.  

Using $x_\mu=-X_\mu/X_x$, we obtain
$$\frac{\partial \dot{X}}{\partial \mu}=\left(X_{x\mu}-\frac{X_{xx}X_\mu}{X_x}\right)v-v_xX_\mu+X_xv_\mu\,.$$
Take $X(x,0)=x$ at $\mu=0$, and for all small $\mu$ keep $X(x,\mu)=x$ in a small interval $|x|\le\delta$ so as to preserve \eqref{eq:dotxmux2}.
First  make $\frac{\partial \dot{X}}{\partial \mu}\ge c'$ for some $c'>c$ at $\mu=0$.  At $\mu=0$
$$\frac{ \partial \dot{X}}{\partial\mu}=X_{x\mu}v-v_xX_\mu+v_\mu=\left(\frac{X_\mu}{v}\right)_xv^2+v_\mu\,.$$
Choose a smooth function $g(x)\ge c'$, and $g(x)=1$ in $|x|\le\delta$.
Let 
$$X=x+\mu v(x,0)\int_0^x\frac{g(x')-v_\mu(x',0)}{v^2(x',0)}\, dx'\,.$$
Then $$\frac{\partial \dot{X}}{\partial \mu}=g\ge c'\,.$$
To make sure that $\frac{\partial \dot{X}}{\partial\mu}\ge c$ for small $\mu$, we make the additional constraint on $g$ that $\int_0^L\frac{g-v_\mu}{v^2}\,dx=0$.

\section*{Appendix C: $x_2$-nullclines for $\mu_2\le-2\delta$}

We show that for $\mu_2\le-2\delta$ the $x_2$ nullcline consists of two $C^1$ graphs $x_2^\pm(x_1)$ in $\frac23\sqrt{-\mu_2}\le\pm x_2\le2\sqrt{-\mu_2}$, of small slope.

Let us treat the case $x_2>0$.  For $0\le x_2\le\frac23\sqrt{-\mu_2}$ we have 
\begin{equation*}
\tilde{v}_2\le-|\mu_2|+\frac49|\mu_2|+\frac{|\mu_2|}{2}+\frac{8K}{27}|\mu_2|^{3/2}=-\frac1{18}|\mu_2|+\frac{8K}{27}|\mu_2|^{\frac32}<0
\end{equation*}
for $|\mu_2|<\left(\frac{3}{16K}\right)^2$ (true, because we already assumed $|\mu_2|<\frac{1}{16K^2}$).  We already showed that $\tilde{v}_2>0$ for $x_2\ge2\sqrt{-\mu_2}$.

For $\frac23\sqrt{-\mu_2}\le x_2\le2\sqrt{-\mu_2}$, using $'$ to denote $\frac{\partial}{\partial x_2}$ , we have
\begin{equation*}
\tilde{v}_2'\ge2x_2-\delta-3Kx_2^2
\ge\frac43{\sqrt{-\mu_2}}-\left(\frac12+12K\right)|\mu_2|\ge\sqrt{-\mu_2}
\end{equation*}
for $|\mu_2|\le\frac{4}{9(1+24K)^2}$,
which we will henceforth assume, as we are interested only in small $\mu$.
So for each $x_1$ there is a unique positive $x_2(x_1)$ at which 
\begin{equation}
\tilde{v}_2(x_1,x_2(x_1))=0.\label{eq:tildev2}
\end{equation}
Using the implicit function theorem, for each $\bar{x}_1$, $x_2(\bar{x}_1)$ extends to a $C^1$ solution of \eqref{eq:tildev2} for $x_1$ in a neighbourhood of $\bar{x}_1$, thus, by uniqueness the function $x_2(x_1)$ is $C^1$.

Applying the chain rule to \eqref{eq:tildev2},
\begin{equation*}
\frac{dx_2}{dx_1}=-\frac{1}{\tilde{v}_2'}\frac{\partial \tilde{v}_2}{\partial x_1}
\end{equation*}
so 
\begin{equation*}
\left|\frac{dx_2}{dx_1}\right|\le\frac{\delta}{\sqrt{-\mu_2}}\,.
\end{equation*}

The case $x_2<0$ is similar.

\section*{Appendix D: A region with $C^1$ invariant circles}

As the tangential dynamics expands a lot for $0<x_1\ll L_1/2$ and then contracts a lot for $0<L_1-x_1\ll L_1/2$, the first thing we do is introduce a new horizontal coordinate in which the expansion and contraction are small.

Let $\bar{v}$ be a positive $C^1$ function of $x_1$, close to $v_1$, and write $\Delta v=v_1-\bar{v}$.  Let $y_1=\int_0^{x_1}\frac{dx}{\bar{v}(x)}$ be the new horizontal coordinate.  It has the interpretation of time from $x_1=0$ using the vector field $\bar{v}$.  From the horizontal dynamics  $\dot{x}_1=v_1(x_1)+\mathcal{O}(\delta)$ we obtain
$$\dot{y}_1=\frac{\dot{x}_1}{\bar{v}}=1+\frac{\Delta v+ \mathcal{O}(\delta)}{\bar{v}}.$$
To evaluate horizontal expansion or contraction in $y_1$,
$$\dot{\delta y_1}=\frac{(\Delta v'+\mathcal{O}(\delta))}{\bar{v}}\,\delta x_1+\frac{\mathcal{O}(\delta)}{\bar{v}}\,\delta x_2-\frac{\bar{v}'}{\bar{v}^2}(\Delta v+\mathcal{O}(\delta))\,\delta x_1\,.$$
Using $\delta y_1=\frac{\delta x_1}{\bar{v}(x_1)}$ we obtain
\begin{equation}
\dot{\delta y_1}=[\Delta v'+\mathcal{O}(\delta)-\frac{\bar{v}'}{\bar{v}}(\Delta v+\mathcal{O}(\delta))]\,\delta y_1+\frac{\mathcal{O}(\delta)}{\bar{v}}\delta x_2\,.\label{eq:deltay1dot}
\end{equation}

The vertical dynamics $\dot{x}_2=v_2(x_2)+\mathcal{O}(\delta)$ has linearisation
\begin{equation}
\dot{\delta x_2} =[v_2'(x_2)+\mathcal{O}(\delta)]\,\delta x_2+\mathcal{O}(\delta)\,\delta x_1
=[v_2'(x_2)+\mathcal{O}(\delta)]\,\delta x_2+\mathcal{O}(\delta)\bar{v}\,\delta y_1\,.\label{eq: deltax2dot}
\end{equation}
Combining \eqref{eq:deltay1dot} and \eqref{eq: deltax2dot}, the slope $s=\frac{\delta x_2}{\delta y_1}$ of a tangent vector evolves by the Ricatti equation
\begin{eqnarray}
\dot{s}&=&\frac{\delta y_1\dot{\delta x_2}-\dot{\delta y_1}\delta x_2}{\delta y_1^2}\nonumber\\
&=&\left[v'_2(x_2)+\mathcal{O}(\delta)-(\Delta v'+\mathcal{O}(\delta))+\frac{\bar{v}'}{v}(\Delta v+\mathcal{O}(\delta))\right]s+\mathcal{O}(\delta)\bar{v}-\frac{\mathcal{O}(\delta)}{\bar{v}}s^2\,.
\label{eq:sdot}
\end{eqnarray}
Each $\mathcal{O}(\delta)$ term in equations (\ref{eq:deltay1dot}--\ref{eq:sdot}) is bounded by $\delta$, because they represent sizes of components of the perturbation or its first partial derivatives.  

Let us specialise to the annulus $A^+_2$ where $v_2'>0$.  To deduce that $C^+_2$ is a $C^1$ circle, it suffices to find conditions on $(\mu_1,\mu_2,\delta)$ such that a cone $|s|\ge s_0$ is forward invariant  and vectors in it are expanded exponentially.  The value $s_0=\sqrt{|\mu_2|}$ will work.  We show the expansion property first.

For $\mu_2<-C\delta$ for $C>1$ we can take the narrower backward invariant annulus with boundaries $x_2=\sqrt{-\mu_2\pm\delta}$ (neglecting the $\mathcal{O}(x_2^3)$ terms in $v_2$ here and throughout, which make small corrections but complicate the formulae).
In this annulus $v_2'(x_2)\ge 2\sqrt{-\mu_2-\delta}$.

For slopes $|s|\ge s_0=\sqrt{-\mu_2}$, using $\delta y_1=\delta x_2/s$ we obtain from equation \eqref{eq: deltax2dot}
$$\frac{\dot{\delta x_2}}{\delta x_2}\ge 2\sqrt{-\mu_2-\delta}-\delta-\frac{B\delta}{\sqrt{-\mu_2}}$$
where $B=\sup|\bar{v}|$.
By hypothesis, $\delta\le|\mu_2|/C$, so
\begin{equation*}
\frac{\dot{\delta x_2}}{\delta x_2}\ge 2\sqrt{|\mu_2|}\sqrt{1-\frac{1}{C}}-\frac{|\mu_2|}{C}-
\frac{B}{C}\sqrt{|\mu_2|}=\sqrt{|\mu_2|}\left(2\sqrt{1-\frac{1}{C}}-\frac{B}{C}\right)-\frac{|\mu_2|}{C}\,.
\end{equation*}
By increasing $C$ a little above $1+B^2/4$  if necessary, we can make this positive, indeed as close as we like to $2\sqrt{|\mu_2|}$.

To make the cone $|s|\ge s_0$ forward invariant we have to choose $\bar{v}$.  First treat the case $\mu_1\ge|\mu_2|$.  Then we can take $\bar{v}=v_1$.  Since $\Delta v=0$ in this case, and without loss of generality considering $s>0$, from equation \eqref{eq:sdot}
\begin{equation*}
\dot{s}\ge\left[2\sqrt{-\mu_2-\delta}-2\delta-\left|\frac{\bar{v}'}{\bar{v}}\right|\delta\right]s-B\delta-\frac{\delta}{\mu_1}s^2\,.
\end{equation*}
Near $x_1=0$, $\bar{v}\sim\mu_1+x_1^2$ 
so $\frac{\bar{v}'}{\bar{v}}\sim\frac{2x_1}{\mu_1+x_1^2}$ 
which is at most $1/\sqrt{\mu_1}$ (achieved at $x_1=\sqrt {\mu_1}$); outside a larger neighbourhood of $0$, $\frac{\bar{v}'}{\bar{v}}$ is bounded independently of $\mu_1$.  Thus on $s=s_0$
\begin{equation*}
\dot{s}\ge\left(2\sqrt{-\mu_2-\delta}-2\delta-\frac{\delta}{\sqrt{\mu_1}}\right)\sqrt{-\mu_2}-B\delta-\frac{\delta}{\mu_1}|\mu_2|
\end{equation*}
Using $\delta\le|\mu_2|/C$ and $\mu_1\ge|\mu_2|$, we see that at $s=s_0$, 
$$\dot{s}\ge2|\mu_2| \sqrt{1-\frac1C}-\frac2C|\mu_2|^\frac32-
\frac{|\mu_2|}C-\frac{B|\mu_2|}C-\frac{|\mu_2|}C\,.$$
Taking $C$ a little larger than $1+(1+B/2)^2$, this is positive, indeed by increasing $C$ we can make it as close as we like to $2|\mu_2|$.

\begin{figure}[htbp] 
   \centering			
   \includegraphics[width=2.5in]{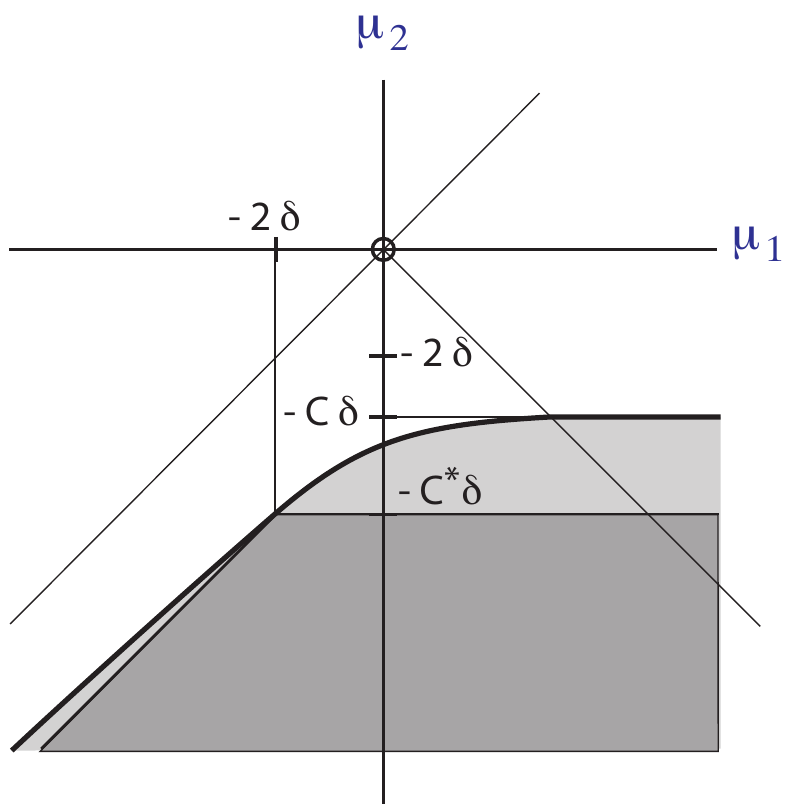} 
   \caption{A region with $C^1$ invariant circles.}
   \label{fig:C1region}
\end{figure}

Next we do the case $\mu_1<|\mu_2|$, where we will find a constraint of the form $\mu_1$ greater than some function of $\mu_2$.  This time, we take $\Delta v$ to be the (negative) constant $\mu_1-|\mu_2|$.  Then $\Delta v'=0$ and $\bar{v}\ge|\mu_2|$.  So (for $s>0$)
\begin{equation*}
\dot{s}\ge\left[2\sqrt{-\mu_2-\delta}-2\delta-\left|\frac{\bar{v}'}{\bar{v}}\right|(|\mu_2|-\mu_1+\delta)\right]s-B\delta-\frac{\delta}{|\mu_2|}\,s^2.
\end{equation*}
Near $x_1=0$, $\bar{v}\sim|\mu_2|+x_1^2$ so $\frac{\bar{v}'}{\bar{v}}\sim\frac{2x_1}{|\mu_2|+x_1^2}$ is at most $\frac1{\sqrt{|\mu_2|}}$; again it is bounded uniformly outside a larger neighbourhood of $0$. So at $s=\sqrt{|\mu_2|}$, 
\begin{equation*}
\dot{s}\ge2\sqrt{-\mu_2-\delta}\sqrt{|\mu_2|}-2\delta\sqrt{|\mu_2|}-(|\mu_2|-\mu_1+\delta)-B\delta-\delta\,.\label{eq:sdot2}
\end{equation*}
This is non-negative if and only if 
\begin{equation}
\mu_1\ge |\mu_2|+(B+2)\delta+2\delta\sqrt{|\mu_2|}-2\sqrt{|\mu_2|-\delta}\sqrt{|\mu_2|}\,.\label{eq:mu12}
\end{equation}
This region is sketched in Figure~\ref{fig:C1region}.
Hence the cone $|s|\ge s_0$ is forward invariant for those $\mu_1$ in $(-|\mu_2|,|\mu_2|)$ satisfying \eqref{eq:mu12}.

The region \eqref{eq:mu12} in which we have obtained a $C^1$ circle has too complicated a formula to carry around with us, but it contains a region of the form
\begin{equation*}
\mu_2\le
\begin{cases}
-C^*\delta\quad\text{for }\mu_1\ge-2\delta\\
-C^*\delta+\mu_1+2\delta\quad\text{for } \mu_1\le-2\delta
\end{cases}
\end{equation*}
as shown in Figure~\ref{fig:C1region}, where $C^*$ is determined by the intersection of \eqref{eq:mu12} with $\mu_1=-2\delta$. It can be written compactly as
$\mu_2\le -C^*\delta$, $\mu_1\ge\mu_2+(C^*-2)\delta$.
We henceforth denote $C^*$ by $C$.

\section*{Appendix E: Corrections to the transit map}    

To estimate the effects on the transit map of corrections to $v_1$ and $v_2$ and of the parameters $\mu_j$ and the perturbation $\delta$, we start from the explicit choice \eqref{eq:v2}
of $v_2(x_2)=\left(\frac{L_2}{\pi}\right)^2\sin^2\frac{\pi x_2}{L_2}$, and take $v_1=\mu_1+x_1^2$ in $|x_1|\le\eta$.  
Then
\begin{eqnarray*}
x_1(t)&=&\frac{x_1(0)}{1-x_1(0)t}\quad\qquad\qquad\qquad\qquad\text{for }\mu_1=0\\
&=&\frac{\sqrt{\mu_1}\tan\sqrt{\mu_1}t+x_1(0)}{1-\frac{x_1(0)}{\sqrt{\mu_1}}\tan\sqrt{\mu_1}t}\qquad\qquad\text{for }\mu_1>0\\
&=&\frac{\sqrt{-\mu_1}\tanh\sqrt{-\mu_1}t+x_1(0)}{1+\frac{x_1(0)}{\sqrt{-\mu_1}}\tanh\sqrt{-\mu_1}t}\qquad\text{for }\mu_1<0\,.
\end{eqnarray*}
Since $t_2\sim2/\eta$ and using $\sqrt{\mu_1}\ll\eta$ and $x_1(0)\in [-\eta,\frac{\eta}{3}]$, we see that the effect of $\mu_1$ on $x_1'=x_1(t_2)$ is $O(\mu_1/\eta)$.

To integrate $\dot{x_2}=\left(\frac{L_2}\pi\right)\sin^2\frac{\pi x_2}{L_2}+\mu_2$ we put $\sigma=-\cot\frac{\pi x_2}{L_2}$ which for $x_2\in[\eta,L_2-\eta]$ increases from $-\sigma_0$ to $+\sigma_0$, where $\sigma_0=\cot\frac{\pi\eta}{L_2}\sim\frac{L_2}{\pi\eta}$.  Then 
$$\dot{\sigma}=\csc^2 \frac{\pi x_2}{L_2}  \cdot \frac{\pi}{L_2}\dot{x_2}=\frac{L_2}{\pi}+\frac{\pi}{L_2}\mu_2(1+\sigma^2)=(\frac{L_2}{\pi}+\frac{\pi\mu_2}{L_2})+\frac{\pi\mu_2}{L_2}\sigma^2\,.$$
So
\begin{eqnarray*}
t_2&=&\frac{2\pi}{L_2}\sigma_0\;\qquad\qquad\qquad\qquad\qquad\qquad\qquad\qquad\text{if }\mu_2=0\\
&=&\frac{2}{\sqrt{(\frac{L_2}{\pi}+\frac{\pi\mu_2}{L_2})\frac{\pi\mu_2}{L_2}}}
\tan^{-1}
\sqrt{\frac{\frac{\pi\mu_2}{L_2}}  {\frac{L_2}{\pi}  +\frac{\pi\mu_2}{L_2}}     }\,
\sigma_0
\qquad\text{if }\mu_2>0\\
&=&\frac{2}{\sqrt{(\frac{L_2}{\pi}+\frac{\pi\mu_2}{L_2})\frac{\pi|\mu_2|}{L_2}}}
\tanh^{-1}
\sqrt{\frac{\frac{\pi|\mu_2|}{L_2}}  {\frac{L_2}{\pi}  +\frac{\pi\mu_2}{L_2}}     }\,
\sigma_0
\quad\text{if }\mu_2<0\,.
\end{eqnarray*}
The dominant change to $t_2$ with $\mu_2$ comes from the Taylor expansions of $\tan^{-1}x \sim x-\frac{x^3}{3}$, $\tanh^{-1}x\sim x+\frac{x^3}{3}$ rather than the $\frac{L_2}{\pi}+\frac{\pi\mu_2}{L_2}$ terms.  So for instance for positive $\mu_2\ll \eta^2$,
\begin{equation*}
t_2\sim\frac{2\pi}{L_2}\cot\frac{\pi\eta}{L_2}-\frac{2}{3}\frac{\pi^3\mu_2}{L_2^3}\cot^3\frac{\pi\eta}{L_2}
\sim\frac{2\pi}{L_2}\sigma_0-\frac23\frac{\pi^3\mu_2}{L_2^3}\sigma_0^3\sim\frac2\eta-\frac23\frac{\mu_2}{\eta^3}\,. 
\end{equation*}
The effect of a change $\Delta t_2$ in $t_2$ on $x_1'$, using the case $\mu_1=0$ and $x_1\in[-\eta,\eta/3]$, is
\begin{equation}\label{eq:Deltax1'}
\Delta x_1'=\frac{x_1}{1-x_1(t_2+\Delta t_2)}-\frac{x_1}{1-x_1t_2}=\frac{x_1^2\Delta t_2}{(1-x_1(t_2+\Delta t_2))(1-x_1 t_2)}=\mathcal{O}(x_1^2\,\Delta t_2).
\end{equation}
So the effect of the change $-\frac{2}{3}\frac{\mu_2}{\eta^3}$ to $t_2$ is $\Delta x_1'=\mathcal{O}\left(\frac{x_1^2\mu_2}{\eta^3}\right)$.

Next we consider the change induced by deforming $v_2$ to an arbitrary $C^3$ vector field with the same second order Taylor expansion and positive away from zero.  We can write such a deformation as $\alpha(x_2)\sin ^3\frac{\pi x_2}{L_2}$ for a bounded function $|\alpha(x_2)|\le\alpha_0$.  Then $t_2$ is given by
\begin{equation*}
t_2=\int_\eta^{L_2-\eta}\frac{dx_2}{v_2+\alpha\sin^3\frac{\pi x_2}{L_2}}\sim\int\frac{dx_2}{v_2}-\int\frac{\alpha\sin^3\frac{\pi x_2}{L_2}}{v_2^2}\,dx_2\,.
\end{equation*}
The leading correction to $t_2$ is 
\begin{eqnarray*}
-\int\frac{\alpha(x_2)\,dx_2}{\left(\frac{L_2}{\pi}\right)^4\sin\frac{\pi x_2}{L_2}}&\le&\left(\frac{\pi}{L_2}\right)^4\alpha_0\left[-\log(\csc\frac{\pi x_2}{L_2}+\cot\frac{\pi x_2}{L_2})\right]^{L_2-\eta}_\eta\nonumber\\
&=&2\left(\frac{\pi}{L_2}\right)^4\log(\sqrt{1+\sigma_0^2}+\sigma_0)\sim 2\left(\frac{\pi}{L_2}\right)^4\alpha_0\log\frac{2L_2}{\pi\eta}\,.
\end{eqnarray*}
Then using \eqref{eq:Deltax1'} the correction $\Delta x_1'=\mathcal{O}(x_1^2\log\frac{1}{\eta})$.

The effect of a perturbation $\delta$ on $\dot{x_1}$ is at most like changing $\mu_1$ by $\delta$.  By our previous analysis this produces $\Delta x_1'=\mathcal{O}\left(\frac{\delta}{\eta}\right)$.  

The effect of a perturbation $\delta$ on $\dot{x_2}$ can be split into two parts. Firstly, there is a part  bounded by $\sin^3\frac{\pi x_2}{L_2}$ which can be absorbed into  $\alpha$
so, assuming the $C^3$ norm of the perturbation is $\mathcal{O}(\delta)$, it changes $\alpha$ by $\mathcal{O}(\delta)$ and  the resulting change to $\Delta x_1'$ is $\mathcal{O}(\delta x_1^2\log\frac{1}{\eta})$.  Secondly, there is a part localised near the box which fits with expansion \eqref{eq:expansion},  thus
$$
\delta_2x_1+\varepsilon_2x_1^2+\mathcal{O}(x_2^3)+\mathcal{O}(\delta|x|^3)\,,
$$
but whose $\mathcal{O}(x_2^3)$ terms can be absorbed into $\alpha$, so this has size at most $\delta_2|x_1|+\delta\eta^3$.  They have an effect like $\mu_2$ so from \eqref{eq:Deltax1'} they produce $\Delta x_1'=\mathcal{O}(|x_1|^2(\delta_2\frac{|x_1|}{\eta^3}+\delta))$.

Lastly the $\mathcal{O}(x_1^3)$ terms in $\dot{x_1}$ produce $\Delta x_1'=\mathcal{O}(\frac{x_1^3}{\eta})$, because $x_1$ does not change much during the transit.

The effects of the different types of change are close to additive, so putting all this together, and using $|\mu_j|\le C\delta$ we obtain that the error in the transit map for $\mu_1=\mu_2=\delta=0$ is $\mathcal{O}(x_1^2\log\frac1\eta)$ and the effects of $\mu_j$ and $\delta$ on the transit map are $\mathcal{O}(\delta/\eta)$.  

Using the implicit function theorem to determine $t_2$ as a function of $x_1(0)$, one can also deduce that the corrections are small in $C^1$.

\section*{Appendix F: Persistence of invariant circles for $\mu_1\le-K\delta^{4/3}$ }  
In this appendix we sketch the proof that there is $K$ large such that for $\mu_1\le-K\delta^{4/3}$ between the outer sne and lower cusped sne curves, the two vertical circles $C^\pm_1$ persist.

Let us treat the region inside the cusp first.

\begin{figure}[h] 
   \centering
   \includegraphics[width=2in]{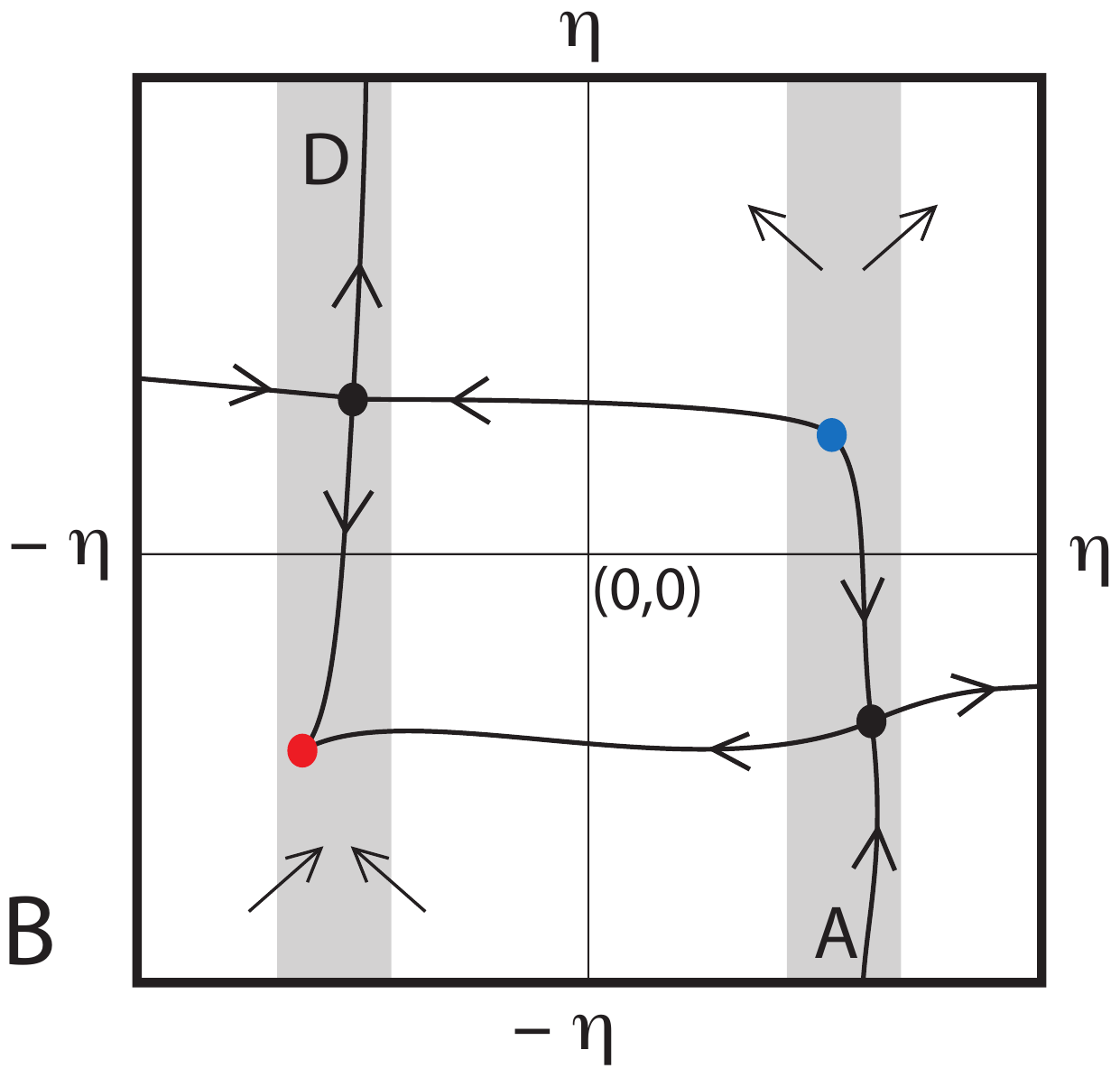} 
   \caption{Box $B$ with strips $\sqrt{|\mu_1|-\delta\eta-\epsilon}\le\mp x_1\le\sqrt{|\mu_1|+\delta\eta+\epsilon}$ for parameters inside the cusp with $|\mu_1|\gg\delta\eta$, $\mu_1<0$.}
   \label{fig:boxAD}
\end{figure}

If $\delta\eta\ll|\mu_1|$, the flow is inward, respectively outward, across the boundaries of the strips 
$$\sqrt{|\mu_1|-\delta\eta-\epsilon}\le\mp x_1\le\sqrt{|\mu_1|+\delta\eta+\epsilon}$$
 in the box $B$ (see Figure~\ref{fig:boxAD}), where $\epsilon$ is a small amount to take care of the remainder terms in $\dot{x}_1$.  
The equilibria all lie in these strips.  The upwards unstable manifold $D$ of the top left saddle and the downward stable manifold $A$ of the lower right saddle stay in the corresponding strips until they exit the box (the notation for branches of saddle manifold is introduced in generality in Figure~\ref{fig:ABCDlabels}).  The transit map moves $D$ by $\mathcal{O}(\delta/\eta)$ so if this is much less than $\sqrt{|\mu_1|}$  then $D$ arrives at $x_2=L_2-\eta$ to the left of $A$.  Because of the way the equilibria connect in the box (Figure \ref{fig:PhasePortaitsBoxDelta2positif}) this implies that $D$ goes to the sink and $A$ to the source, thus forming our two vertical invariant circles.  The two conditions $\delta\eta\ll|\mu_1|$, $\frac{\delta}{\eta}\ll\sqrt{|\mu_1|}$ can be achieved for the greatest region of $|\mu_1|$ if we choose $\eta\sim\delta^{1/3}$ (a choice we have already found useful) and then the conditions hold for $|\mu_1|\ge K\delta^{4/3}$ for large enough $K$ as claimed.

\begin{figure}[h] 
   \centering
   \includegraphics[width=2in]{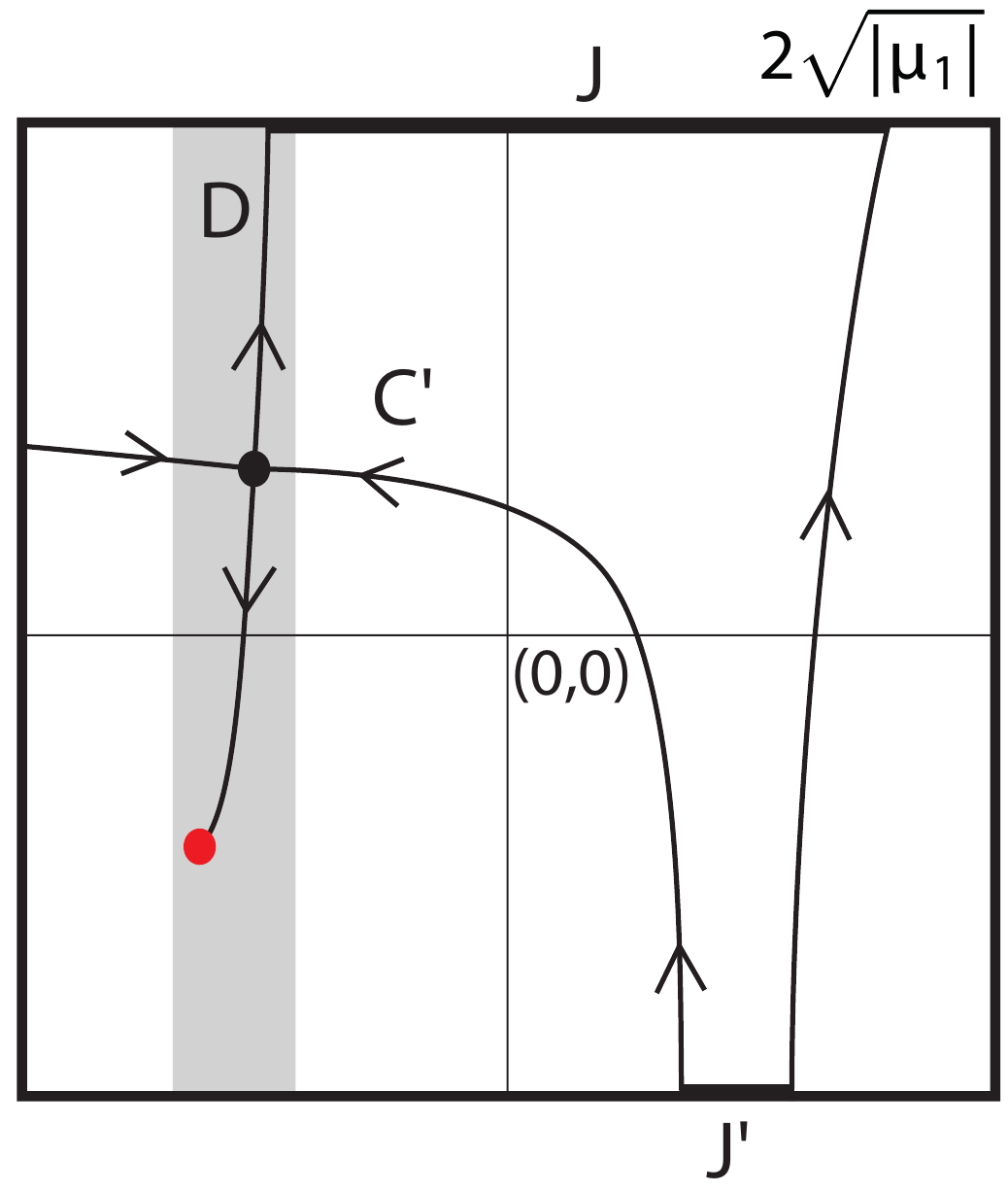} 
   \caption{Box $B$ with strip $-\sqrt{|\mu_1|+\delta\eta+\epsilon}\le x_1\le-\sqrt{|\mu_1|-\delta\eta-\epsilon}$ and the backwards orbit of $(2\sqrt{|\mu_1|},\eta)$, for parameters between the outer sne curve and the upper cusped sne curve with $|\mu_1|\gg\delta\eta$, $\mu_1<0$.}
   \label{fig:boxDC'}
\end{figure}

Next we treat the region between the outer sne and the upper cusped sne curves.  The lefthand strip contains a saddle and a sink and the saddle manifold $D$ reaches $x_2=L_2-\eta$ to the left of $-\sqrt{|\mu_1|-\delta\eta-\epsilon}+\mathcal{O}(\delta/\eta)$.  We now need to bound where the rightward branch $C'$ of stable manifold of the saddle goes.  We claim that for $|\mu_1|\gg\delta^2$ it crosses $x_1=0$ inside the box and hence exits $x_2=-\eta$ with $x_1>0$ as in Figure~\ref{fig:boxDC'}.

To sketch why this is so, make the approximations
\begin{equation*}
\begin{cases}
\dot{x}_1&=\mu_1+x_1^2\\
\dot{x}_2&=\delta_2(x_1+\sqrt{|\mu_1|})
\end{cases}
\end{equation*}
near the saddle ($x_1\approx-\sqrt{|\mu_1|}$, $x_2\approx+\sqrt{\delta_2\sqrt{|\mu_1|}-\mu_2}$).
The first has solution $x_1(t)=-\sqrt{|\mu_1|}\tanh\sqrt{|\mu_1|}t$ with origin of time chosen to correspond to $x_1(0)=0$.  Then the second equation is $\dot{x}_2=\delta_2\sqrt{|\mu_1|}(1-\tanh\sqrt{|\mu_1|}t)$ which has solution $x_2(t)=x_2(+\infty)-\delta_2\log(1+e^{-2\sqrt{|\mu_1|}t})$.  Thus $x_2(0)\approx\sqrt{\delta_2\sqrt{|\mu_1|}-\mu_2}-\delta_2\log2$ which is well inside the box for $\delta_2\ll\eta$.  The principal correction to $\dot{x}_1$ is $\delta_1x_2$ but this is by the hypothesis $\delta\eta\ll|\mu_1|$ small compared to $\mu_1$.  The principal correction to $\dot{x}_2$ is $x_2^2+\mu_2-\delta_2\sqrt{|\mu_1|}$, but for $t$ in the interval $[0,+\infty)$ this makes a negligible change to $x_2$.

Thus for $\delta/\eta\ll\sqrt{|\mu_1|}$, $D$ arrives  to the left of $C'$ and hence goes into the sink, making an invariant circle $C^-_1$.

To complete the analysis of this case we construct a periodic orbit for $C^+_1$.  The interval $J$ on $x_2=+\eta$ between $D$ and $2\sqrt{|\mu_1|}$ under the backward flow exits $x_2=-\eta$ by an interval $J'$ between $C'$ and something a little to the right of $+\sqrt{|\mu_1|}$.  This is because 
$\dot{x}_1=\mu_1+x_1^2$
in backwards time contracts $x_1>\sqrt{|\mu_1|}$ towards $\sqrt{|\mu_1|}$ roughly like $e^{\sqrt{|\mu_1|}t}$.  
To estimate the time $t$ it takes $x_2$ to flow backwards from $+\eta$ to $-\eta$, use $\dot{x}_2\approx\mu_2+x_2^2+\delta_2x_1$ and note that $|\mu_2|\lesssim\delta_2\sqrt{|\mu_1|}$ between the sne curves and $x_1\le2\sqrt{|\mu_1|}$, so 
$|\dot{x}_2|\le3\delta_2\sqrt{|\mu_1|}+x_2^2$.  The solution of $\dot{x}=a^2+x^2$ is $x(t)=a\tan at$ so takes time $\sim\frac{\pi}{2a}$ to cross $|x|\lesssim a$.  Hence $x_2$ takes at least time $\frac{1}{\sqrt{3\delta\sqrt{|\mu_1|}}}$.
This time greatly exceeds $1/\sqrt{|\mu_1|}$.

The backwards transit map from $x_2=-\eta$ to $x_2=-L_2+\eta$ moves points of $J'$ by $\mathcal{O}(\delta/\eta)$ so if this is less than $\sqrt{|\mu_1|}$  the interval $J$ is mapped by the backwards flow strictly inside itself and so has at least one attracting fixed point, making a periodic orbit.  With some estimates of its derivative we could establish uniqueness.

\section*{Appendix G: Unfolding of an $E$ point}
\begin{figure}[h] 
   \centering
   \includegraphics[width=4in]{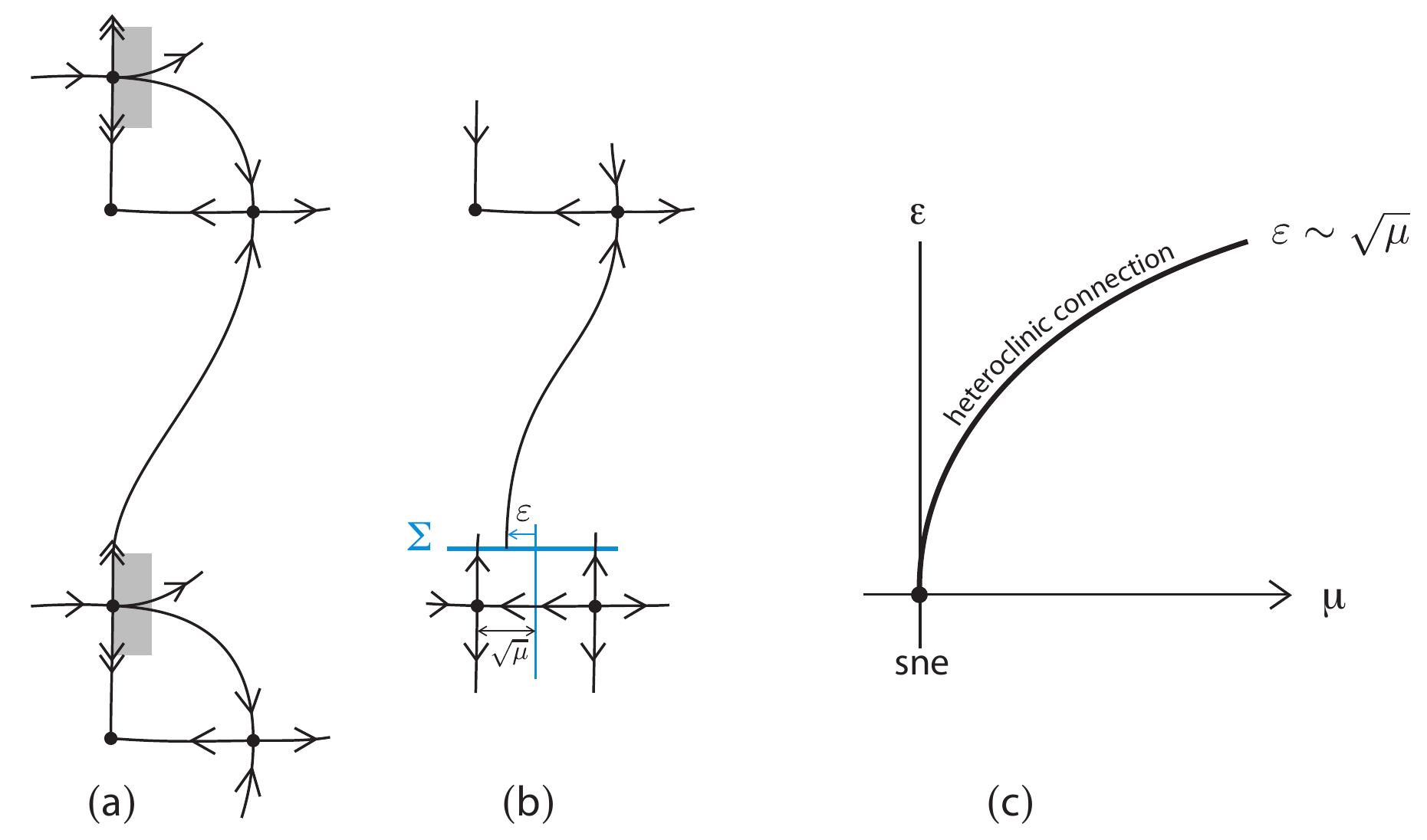} 
   \caption{(a) Phase portrait for an $E$ point; (b) Effects of unfolding parameters $\mu,\eps$; (c) Bifurcation diagram.}
   \label{fig:Ept}
\end{figure}

An $E$ point is the codimension-two situation with an elementary saddle-node equilibrium whose strong unstable manifold connects to the stable manifold of a saddle whose other branch of stable manifold lies in the repelling half-plane of the saddle-node, forming a homotopically non-trivial cycle (or the time-reverse of this situation).  See Fig.~\ref{fig:Ept}(a).

It can be unfolded by a parameter $\mu$ to unfold the saddle-node equilibrium and a parameter $\eps$ to displace the intersection of the stable manifold of the saddle with a transverse section $\Sigma$.  See Fig.~\ref{fig:Ept}(b).  The resulting bifurcation diagram is shown in Fig.~\ref{fig:Ept}(c).


\begin{thebibliography}{XXXXX}
\bibitem[AP]{AP} Andronov AA, Pontryagin LS, Coarse systems (in Russian), Doklady Akademii Nauk SSSR 14 (5) (1937) 247--250
\bibitem[BGKM]{BGKM} Baesens C, Guckenheimer J, Kim S, MacKay RS, Three coupled oscillators: mode-locking, global bifurcations and toroidal chaos, Physica D 49 (1991) 387--475
\bibitem[BM]{BM} Baesens C, MacKay RS, Resonances for weak coupling of the unfolding of a saddle-node periodic orbit with an oscillator, Nonlinearity 20 (2007) 1283--1298
\bibitem[DRS]{DRS} Dumortier F, Roussarie R and Sotomayor J, Generic 3-parameter families of vector fields on the plane, unfolding a singularity with nilpotent linear part. The cusp case of codimension 3, Ergod. Theory Dyn. Syst. 7 (1987)
375--413.
\bibitem[EK]{EK} Ermentrout GB, Kopell N, Parabolic bursting in an excitable system coupled with a slow oscillation, SIAM J Appl Math 46 (1986) 233--253
\bibitem[F]{F} Fenichel N, Persistence and Smoothness of Invariant Manifolds for Flows, Indiana Univ. Math. J. 21 (3) (1972) 193Ð226
\bibitem[GG]{GG} Golubitsky M, Guillemin V, Stable Mappings and Their Singularities, Graduate Texts in Mathematics, Vol. 14, Springer-Verlag (Berlin, 1973)
\bibitem[GK]{GK} Guckenheimer J, Khibnik I, Torus maps from weak coupling of strong resonances. Methods of qualitative theory of differential equations and related topics, Supplement, 205--218, Amer. Math. Soc. Transl. Ser. 2, 200, Amer. Math. Soc., Providence, (RI, 2000).
\bibitem[HI]{HI} Hoppensteadt FC, Izhikevich EM, Weakly connected neural networks (Springer, 1997)
\bibitem[Iz]{Iz} Izhikevich EM, Dynamical Systems in Neuroscience, MIT Press (2007)
\bibitem[LHM]{LHM} Levi M, Hoppensteadt FC, Miranker WL, Dynamics of the Josephson Junction, Q. Appl. Math 36 (1978) 167--198
\bibitem[M]{M} Mazo JJ, Localized excitations in Josephson arrays. Part I: Theory and modeling, in Energy localisation and transfer, eds Dauxois T, Litvak-Hinenzon A, MacKay R, Spanoudaki A (World Sci, 2004), 193--245
\bibitem[NMW]{NMW} Nicholls JG, Martin AR, Wallace BG, From Neuron to Brain, Sinauer (1992)
\bibitem[RE]{RE} Rinzel J, Ermentrout, Analysis of neural excitability and oscillations, in: Methods in neuronal modelling: from ions to networks, 2nd edition, Koch C and Segev I eds, MIT Press (1998)
\bibitem[Sch]{Sch} Schechter S, The saddle-node separatrix-loop bifurcation, SIAM J Math Anal 18 (1987) 1142--56
\bibitem[U]{U} Ustinov AV, Localized Excitations in Josephson Arrays. Part II: Experiments, in Energy localisation and transfer, eds Dauxois T, Litvak-Hinenzon A, MacKay R, Spanoudaki A (World Sci, 2004), 247--271
\end{thebibliography}
\end{document}